\documentclass[11pt]{amsart}
\usepackage{amsmath, amssymb, amsfonts, amsthm, fullpage}
\usepackage{graphics,graphicx}
\usepackage{accents}
\usepackage{color}

\newcommand{\R}{\mathbb{R}}

\newcommand{\into}{\int_{\Omega}}
\newcommand{\intg}{\int_{\partial\Omega}}
\newcommand{\Hoo}{\accentset{\scriptstyle o}{H}^1}
\newcommand{\Hd}[1]{H(\text{div},#1)}

\newcommand{\Hdoo}[1]{\accentset{\scriptstyle o}{H}(\text{div},#1)}

\newcommand{\Voo}{\accentset{\scriptstyle o}V}

\newcommand{\bx}{\mathbf{x}}
\newcommand{\bn}{\mathbf{n}}
\newcommand{\bu}{\mathbf{u}}
\newcommand{\bw}{\mathbf{w}}
\newcommand{\bv}{\mathbf{v}}
\newcommand{\bF}{\mathbf{f}}

\newcommand{\bphi}{\boldsymbol{\phi}}
\newcommand{\dx}{\text{d}\mathbf{x}}
\newcommand{\ds}{\text{d}\mathbf{s}}
\newcommand{\dt}{\Delta t}
\newcommand{\pd}[2]{\partial_{#1}{#2}}

\newcommand{\Ac}{\mathcal{A}}
\newcommand{\Bc}{\mathcal{B}}
\newcommand{\Cc}{\mathcal{C}}
\newcommand{\Ec}{\mathcal{E}}
\newcommand{\Nc}{\mathcal{N}}
\newcommand{\Pc}{\mathcal{P}}

\newcommand{\Th}[1]{\mathcal{T}_{#1}}
\newcommand{\dOmega}{\partial\Omega}
\newcommand{\Hdiv}{\text{H(div)}}
\newcommand{\Div}{\nabla\!\cdot\!}
\newcommand{\Curl}{\text{curl}}
\newcommand{\la}{\langle}
\newcommand{\ra}{\rangle}
\newcommand{\tbn}[1]{{\left\vert\kern-0.25ex\left\vert\kern-0.25ex\left\vert #1 \right\vert\kern-0.25ex\right\vert\kern-0.25ex\right\vert}}

\newtheorem{remark}{Remark}[section]

\AtBeginDocument{%
  \setlength{\belowdisplayskip}{5pt} \setlength{\belowdisplayshortskip}{5pt}
  \setlength{\abovedisplayskip}{5pt} \setlength{\abovedisplayshortskip}{5pt}
}

\begin{document}

\title{Boussinesq-Peregrine water wave models and their numerical approximation}

\author[Th. Katsaounis]{Theodoros Katsaounis}
\address{CEMSE, KAUST, Thuwal, Saudi Arabia \\ \& Dept. of Math. and Applied Mathematics, Univ. of Crete, Heraklion, Greece \\ \& IACM, FORTH, Heraklion, Greece}
\email{theodoros.katsaounis@kaust.edu.sa}

\author[D. Mitsotakis]{Dimitrios Mitsotakis}
\address{School of Mathematics and Statistics, Victoria University of Wellington, Wellington, New Zealand }
\email{dimitrios.mitsotakis@vuw.ac.nz}

\author[G. Sadaka]{Georges Sadaka}
\address{Laboratoire de Math\'ematiques Rapha\"el Salem, Universit\'e de Rouen Normandie, CNRS UMR 6085,
Avenue de l'Universit\'e, BP 12, F-76801 Saint-\'Etienne-du-Rouvray, France}
\email{georges.sadaka@u-picardie.fr}

\begin{abstract}
In this paper we consider the numerical solution of Boussinesq-Peregrine type systems by the application of the Galerkin finite element method. 
The structure of the Boussinesq systems is explained and certain alternative nonlinear and dispersive terms are compared. A detailed study of the convergence properties of the standard Galerkin method, using various finite element spaces on unstructured triangular grids, is presented.  
 Along with the study of the Peregrine system, a new Boussinesq system of BBM-BBM type  is derived. The new system has the same structure in its momentum equation but differs slightly in the mass conservation equation compared to the Peregrine system. Further, the finite element method applied to the new system has better convergence properties, when used for its numerical approximation. 
 Due to the lack of analytical formulas for solitary wave solutions for the systems under consideration, a Galerkin finite element method combined with the Petviashvili iteration is proposed for the numerical generation of accurate approximations of line solitary waves. Various numerical experiments related to the propagation of solitary and periodic waves over variable bottom topography and their interaction with the boundaries of the domains are presented.  We conclude that both systems have similar accuracy when approximate long waves of small amplitude while the Galerkin finite element method is more effective when applied to BBM-BBM type systems.
\end{abstract}
\maketitle
%
%
%
%
%
%
\section{Introduction}
Two identifying physical properties of water waves are the nonlinearity and dispersion. Dispersion refers to the property of waves of different wavelength to travel with different phase speed, while the nonlinearity is related with the way waves interact with each other and also with the fact that waves of larger amplitude propagate faster than those of smaller amplitude. The balanced combination of these two properties is the basic ingredient for the existence of solitary waves. In fluid mechanics the propagation of water waves is described by the Euler equations of ideal fluids (inviscid and irrotational) \cite{Whitham2011}. Due to the complexity of the Euler equations we usually rely on numerical simulations of simplified mathematical models consisting of nonlinear and dispersive partial differential equations when studying the properties of water waves. Depending on the nature of the waves we study,  there are different mathematical models to describe their propagation. For example there are different mathematical models for long waves, short waves, small amplitude or large amplitude waves, waves in small, intermediate or large water depth.

\begin{figure}[ht!]
  \centering
  \includegraphics[width=0.8\columnwidth]{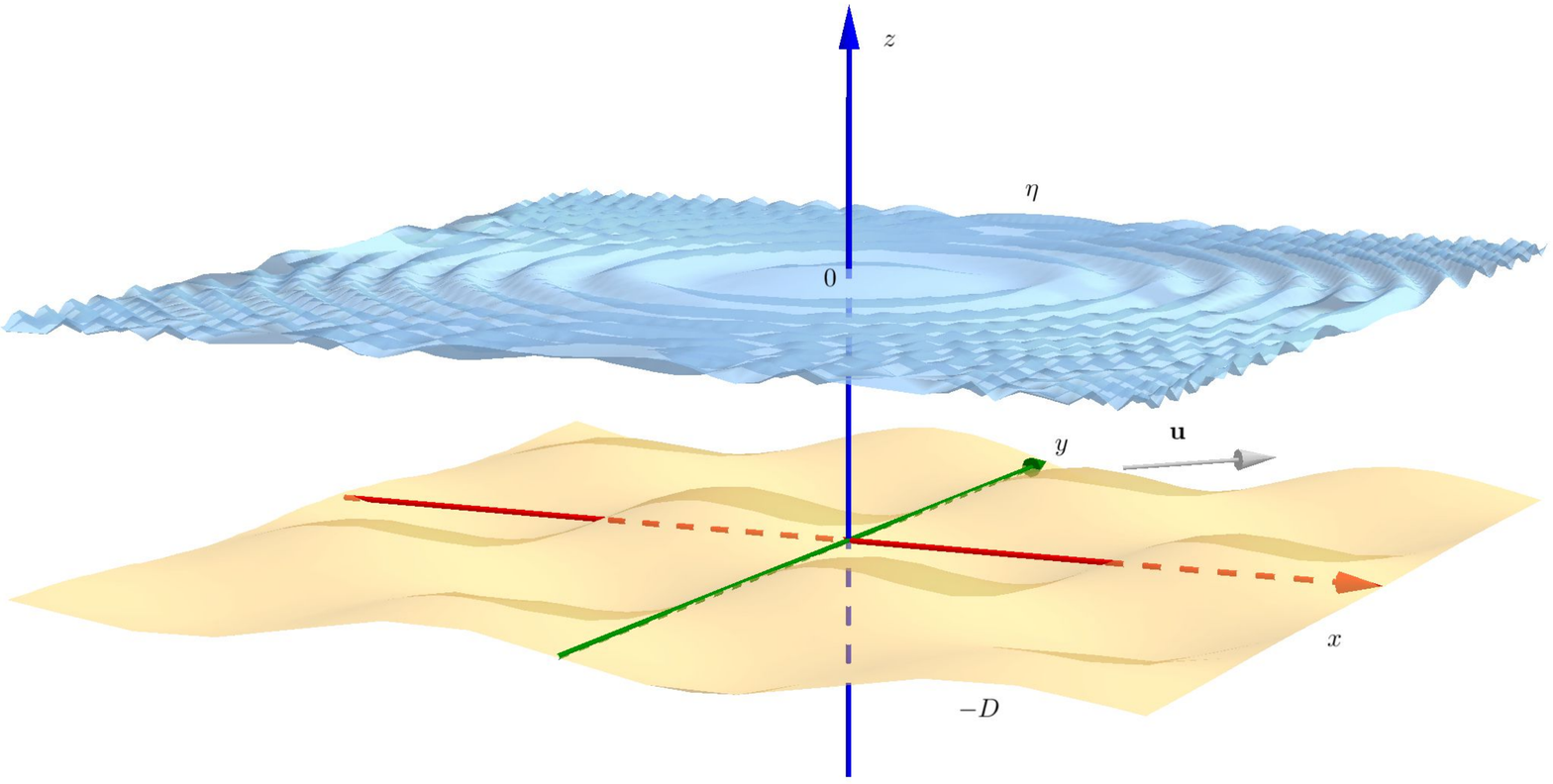}
  \caption{Sketch of the physical domain}
  \label{fig:domain}
\end{figure}

In this paper we focus on some mathematical models derived especially to describe the propagation of small amplitude, long waves compared to the depth of the water. If $A$ denotes a typical wave amplitude in a general state, $D_0$ a typical water depth and $\lambda$ a typical wavelength, then waves that satisfy $\varepsilon:=A/D_0\ll 1$ and $\sigma:=D_0/\lambda\ll 1$ are characterized as long waves of small amplitude in the shallow water regime. Moreover, usually it is required that the Stokes number is $S:=\varepsilon/\sigma^2=O(1)$. Examples of such waves are the important cases of solitary waves and tsumamis in deep ocean waters. The mathematical modelling of water waves of small amplitude in a realistic environment with general bathymetry was initiated by Peregrine \cite{Per67} by deriving asymptotically the so called Peregrine system 
\begin{equation}
\label{eq:Peregrin}
\begin{aligned}
&\eta_t+\nabla\cdot[(D+\eta)\bu]=0\ , \\
&{\bf u}_t+g\nabla\eta+(\bu\cdot \nabla\cdot)\bu -\frac{1}{2}D\nabla(\nabla\cdot(D\bu_t))+\frac{1}{6}D^2\nabla(\nabla\cdot\bu_t)=0\ , 
\end{aligned}
\end{equation}
where $\eta(\bx,t)=\eta(x,y,t)$  denotes the free-surface elevation above an undisturbed level and $\bu(\bx,t)=(u(x,y,t),v(x,y,t))$ is the depth averaged horizontal velocity vector field of the water over a variable bottom  with depth $D(\bx)=D(x,y)>0$. A sketch of the domain is depicted in Figure \ref{fig:domain}. 
It is noted that due to the mean value theorem of integral calculus, $\bu$ can be thought of as the horizontal velocity of the fluid at a certain height above the bottom.

Peregrine's system is known to possess solitary wave solutions similar to the Euler equations \cite{DM2008} and served as a base for the study of water waves, \cite{Debnath2011}. As a mathematical model,  it can be derived using asymptotic techniques and thus has a different linear dispersion relation from that of the Euler equations. Although it is expected to model accurately long waves of small amplitude, there are special cases where the linear dispersion relation plays important role for the accurate approximation of water waves, \cite{Dingemans1994}. Based on the idea of evaluating the horizontal velocity of the fluid at any depth \cite{BS1976}, an alternative Boussinesq system was derived  in the same regime with improved dispersion relation characteristics \cite{Nwogu93} 
\begin{equation}
\label{eq:nwogu}
\begin{aligned}
& \eta_t+\nabla\cdot((D+\eta)\bu)+\nabla\cdot [\tilde{a}D^2 \nabla(\nabla\cdot(D\bu))+\tilde{b}D^3\nabla(\nabla\cdot\bu)]=0\ , \\
& {\bf u}_t+g\nabla\eta+ (\bu\cdot\nabla)\bu+\tilde{c}D\nabla(\nabla\cdot(D\bu_t))+\tilde{d}D^2\nabla(\nabla\cdot\bu_t)=0\ , 
\end{aligned}
\end{equation}
where now $\bu$ denotes the horizontal velocity vector field evaluated at height $z=-D+\theta(\eta+D)$ for $\theta\in[0,1]$, and $\tilde{a}=\theta-1/2$, $\tilde{b}=1/2[(\theta-1)^2-1/3])$, $\tilde{c}=\theta-1$ and $\tilde{d}=1/2(\theta-1)^2$ (see also \cite{mits2009,Nwogu93}).  

The finite element method for the numerical solution of the Peregrine system \eqref{eq:Peregrin} has been introduced in several occasions to study practical problems in coastal engineering \cite{Antunes1993,AQ1998,LLYL1999,SSS2004,ZW2008}. Due to the simplicity of the mathematical formulation of  Peregrine's system, the formulation of the standard Galerkin finite element method in a bounded domain and with slip-wall boundary conditions such as $\bu\cdot\bn=0$,  with $\bn$ being the unit outward normal vector on the boundary of the computational domain, is straightforward. On the other hand its convergence properties for the specific model still remain unknown. 

The standard Galerkin finite element method for the Peregrine system has been studied analytically only in one-dimension and has been proven that its convergence is suboptimal when the same finite element space is used for both dependent variables and it is optimal when different finite element spaces for the approximation of the free-surface elevation and the horizontal velocity are considered  \cite{AD2013, AD2010, AD2012}. The two-dimensional case is different though, and the main reason is that the regularization operator in the momentum conservation of Peregrine's system is the $I-\nabla(\nabla \cdot )$. This operator in one-dimension is identical to the Laplace operator, but in two dimensions is not fully-elliptic as some of the second derivatives are missing. For example, if the solution $\bu$ is restricted in the $\Curl$ vector field space,  then this operator is just the identity.  For this reason we will call this operator incomplete-elliptic operator. Its inversion then can create implications as far as it concerns the regularity of the solution as well as the convergence properties of the numerical method.  It is noted that Nwogu's system (with optimal in some sense linear dispersion properties \cite{Nwogu93}) contains the same incomplete-elliptic operator in the momentum conservation equation and thus similar behaviour in the convergence properties of the variable $\bu$ is expected. 

Another reason for observing suboptimal convergence rates, in the case of Peregrine's system is the form of the mass conservation equation \eqref{eq:Peregrin}, which resembles a first order hyperbolic equation. While in Nwogu's system the mass conservation equation is not a hyperbolic equation,  the presence of the third-order spatial derivative is expected to cause different problems in the finite element semidiscretization and especially complications on the implementation of additional boundary conditions \cite{WK1995}. For the numerical implementation of a finite element method for Nwogu's system we refer to \cite{WB2002}. These peculiarities have been realised in \cite{mits2009} where Boussinesq systems were modified by using the irrotationality condition that governs ideal fluids with potential flow. A particular regularized system of BBM-BBM type was derived in \cite{mits2009} in the form
\begin{equation}
\label{eq:BBMold}
\begin{aligned}
&\eta_t+\nabla\cdot[(D+\eta)\bu]-\frac{1}{6}\nabla\cdot(D^2\nabla\eta_t)=0\ , \\
&\bu_t+g\nabla\eta+(\bu\cdot \nabla)\bu-\frac{1}{6}D^2\Delta\bu_t=0\ .
\end{aligned}
\end{equation}
This system is analogous to the BBM-BBM system derived in \cite{BCL2015} in the flat bottom case. Although the finite element method for the specific system appears to converge with optimal rate, the disadvantage of the specific regularized system is that the inversion of the elliptic operator requires non-slip wall boundary conditions of the form $\bu=0$ or the computation of the tangential component of the velocity vector field on the boundary of the domain  \cite{DMS2007,dms2009}. Non-slip wall boundary conditions can be restrictive and on the other hand the computation of the tangential component of the velocity on the boundary of the domain can be rather complicated. On the other hand, such BBM-BBM type system was shown to have very good performance in studies of water waves over variable bottom \cite{Senth2016,HK2019,mits2009} and for this reason certain improvements can be made so as to make slip-wall boundary conditions easily applicable.

In this paper we study computationally the convergence properties of the standard Galerkin finite element method for the Peregrine system \eqref{eq:Peregrin} and we show that the effects of the hyperbolic form of the mass conservation results in suboptimal convergence rates for the surface elevation. Optimal convergence rates can be obtained only by using different finite element spaces between the free-surface elevation and the velocity field. Furthermore, we show by numerical means that the incomplete-elliptic  operator $I-\nabla(\nabla \cdot )$ results also in suboptimal convergence rates in the $L^2$-norm and optimal in $\Hdiv$-norm. We also show that optimal convergence rates for the velocity variable can be achieved by replacing the incomplete-elliptic regularization operator in the momentum conservation equation by the Laplace operator as in system \eqref{eq:BBMold}. This substitution can be explained using asymptotic arguments. Of course such modification will inflict again restrictive zero Dirichlet boundary condition on the velocity vector field. We also show that the nonlinearity has no implication on the convergence properties of the numerical method. In particular, modifying  the nonlinearity that appears in the momentum conservation equation to one given in conservative form using the same asymptotic arguments we observed the same convergence properties.

Finally, we derive a new, appropriate for slip-wall boundary conditions, form of the BBM-BBM system, \eqref{eq:BBMold}. The new system appears to have optimal convergence rates for the free-surface elevation variable. The momentum conservation equation retains the same mathematical formulation as in the Peregrine system and thus  slip-wall boundary conditions can be used in a straightforward manner. The new system appears to be highly accurate in the appropriate shallow water regime and due to the better convergence properties it can be used as an alternative system for the computations of long water waves of small amplitude. The numerical models are validated against standard benchmarks where the numerical solutions are compared with laboratory data. For the numerical generation of accurate approximation of solitary waves in the experiments we describe a Galerkin finite element method combined with the Petviashvili iteration \cite{Petv1976} for the solution of the resulting nonlinear equations. 

The structure of the paper is as follows: In Section \ref{sec:models} after reviewing the Peregrine system and its properties, we derive the version of the Peregrine system with the standard elliptic regularization operator and the modified nonlinearity. Special attention is being paid to the incomplete-elliptic operator and its numerical solution in Section \ref{sec:feapprox}. In the same section we study the convergence properties of the various Boussinesq systems. The numerical generation of solitary waves is presented in Section \ref{sec:solitwaves} and further testing of the accuracy of the numerical methods on the propagation of solitary waves are considered. Section \ref{sec:numerexp} presents standard benchmark experiments for Boussinesq systems showing that all models with the same asymptotic behaviour result to very similar physically sound solutions. Conclusions and further developments are discussed in Section \ref{sec:conclusions}.
%
%
%
%
%
%
%
%
\section{Mathematical Models}\label{sec:models}
\subsection{Notation}
Let $\Omega\subset \R^2$ be a bounded set with smooth boundary $\Gamma=\dOmega$ and
let $L^2(\Omega)$ denotes standard Lebesque space of square integrable functions in $\Omega$. We shall use the standard notation for the Sobolev spaces defined on $\Omega$ and their corresponding norms:
\begin{align*}
& V:=H^1(\Omega) = \left\{\phi\in L^2(\Omega) : \nabla\phi\in L^2(\Omega) \right\}, \quad \|\phi\|_1^2 = \|\phi\|^2 + \|\nabla\phi\|^2 ,  
\end{align*}
where $\|\cdot\|$ denotes the norm in $V$,  $(\cdot \, ,\cdot)$ is the corresponding inner product in $\Omega$ and $\la \cdot\, , \cdot \ra$ the inner product on $\Gamma$. 
The extension to vector fields is straightforward, i.e. if $\bu, \bv\in L^2(\Omega)^2$ then $(\bu,\bv)=\sum_{i=1}^2 (u_i,v_i)$ and $\|\bu\|^2 = \sum_{i=1}^2 \|u_i\|^2 $. 
We shall also use the following subspace of $H^1(\Omega)$:
\begin{equation*}
\Voo:=\Hoo(\Omega) =\left\{\bv\in H^1(\Omega)^2 : \bv\cdot\bn = 0  \right\}\ ,
\end{equation*} 
where $\bn$ denotes the outward unit normal to $\Gamma$.  Further we shall also work with the Sobolev spaces
\begin{align*}
& \Hd{\Omega} = \left\{\bv\in L^2(\Omega)^2 : \Div\bv\in L^2(\Omega)\right\}, \
 \Hdoo{\Omega}= \left\{\bv\in\Hd{\Omega} : \bv\cdot \bn = 0\ \text{on}\ \Gamma \right\}, 
\end{align*}
with norm 
\begin{equation}
\label{hdivnorm}
\|\bv\|_{\Hdiv}^2 = \|\bv\|^2 + \|\Div\bv\|^2\ .
\end{equation}
The bottom topography is described by $D(\bx)$ and is assumed to be independent of time. Further we assume that there are positive constants $D_0, D_m$ such that $0< D_0 \le D(\bx) \le D_m$. 
\subsection{An elliptic operator}
An interesting  special case of \eqref{eq:Peregrin} is when the bottom topography is flat : $D(\bx)=D_0>0$. In this case the system simplifies considerably and reads
\begin{equation}
\label{prgH0}
\begin{aligned}
& \pd{t}{\eta} + \nabla\cdot\left((D_0+\eta)\bu  \right)  = 0\ , \\
&\partial_t\left(\bu - \frac13 D_0^2 \ \nabla\left(\Div\bu \right)\right) + \left(\bu\cdot\nabla \right)\bu + g\nabla \eta  = 0\ .\\
\end{aligned}
\end{equation}
From systems \eqref{eq:Peregrin} and  \eqref{prgH0} it becomes apparent that the incomplete-elliptic operator 
\begin{equation}
\label{EllOp}
\Ec(\bu) = \bu - \frac12 D \nabla\left(\nabla\cdot(D\bu) \right) + \frac16 D^2 \nabla\left(\nabla\cdot\bu \right)\, \quad \text{or}\quad  \Ec(\bu) = \bu -   \frac13 D_0^2 \nabla(\Div\bu)\ , 
\end{equation}
plays an important role in its solvability. Indeed the invertibility of $\Ec$ will guarantee the solution of the corresponding discrete problem as we shall see next. 
For the operator $\Ec$,  we consider the following boundary value problem : 
\begin{equation}
\label{elliptic}
	\begin{aligned}
	\Ec(\bu) & =  \bF \quad \text{ in } \Omega\ , \\
		\bu\cdot\bn & = 0  \quad \text{ on }\Gamma\ , 
	\end{aligned}
\end{equation}
where $\bF\in L^2(\Omega)^2$. A weak formulation for \eqref{EllOp} is : we seek $\bu\in \Voo$ such that 
\begin{equation}
\label{weakelliptic0}
\Ac(\bu,\bv):=(\Ec(\bu), \bv) = (\bF, \bv), \quad \forall\,\bv\in\Voo\ . \\
\end{equation}
Integrating by parts and using the boundary condition we get : 
\begin{equation}
\label{weakelliptic1}
\begin{aligned}
\Ac(\bu,\bv) & =\into \bu\cdot\bv\,\dx + \frac12 \into \Div(D\bu)\, \Div(D\bv)\, \dx  - \frac16 \into \Div\bu\,\Div(D^2\bv)\, \dx \\
& =  \into \bu\cdot\bv\,\dx + \frac13 \into \Div(D\bu)\, \Div(D\bv)\, \dx + \frac16 \into (\bu\cdot\nabla D) (\bv\cdot\nabla D) \, \dx  \\ 
& + \frac16 \into (\bu\cdot\nabla D) \Div(D\bv) \, \dx - \frac16 \into \Div(D\bu)(\bv\cdot\nabla D)\, \dx
\quad \forall\,\bv\in\Voo \ ,
\end{aligned}
\end{equation}
which is an almost symmetric bilinear form. Then 
 \begin{equation}
\label{weakelliptic2}
\begin{aligned}
\Ac(\bu,\bu) & = \into \bu\cdot\bu\,\dx + \frac13 \into (\Div(D\bu))^2\, \dx + \frac16 \into (\bu\cdot\nabla D)^2\,\dx \\
& =  \into D^{-2} (D\bu)\cdot(D\bu)\,\dx + \frac13 \into (\Div(D\bu))^2\, \dx + \frac16 \into (\bu\cdot\nabla D)^2\,\dx \\
& \ge \min(D_m^{-2}, \frac13) ( \|D\bu\|^2 + \|\Div(D\bu)\|^2) =: C_D^e \|D\bu\|_{\text{H(div)}}^2\ ,
\end{aligned}
\end{equation}
which shows that $\Ac(\cdot, \cdot)$ is coercive with respect to the $\|D\cdot\|_{\Hdiv}$ norm. It is also straightforward to show from \eqref{weakelliptic1} that $\Ac(\cdot, \cdot)$ is continuous in the same norm 
\begin{equation*}
| \Ac(\bu, \bv)| \le C_D^c \|D\bu\|_{\Hdiv} \|D\bv\|_{\Hdiv}  \ \text{ where }\  C_D^c = \max(D_0^{-2}, \frac56)\ , 
\end{equation*}
thus there exist a unique weak solution to \eqref{weakelliptic0} which also satisfies  $\|D\bu\|_{\Hdiv} \le C_f \|\bF\|, \quad C_f = (D_0 C_D^e)^{-1}$.
\begin{remark}
If we assume in addition that  $D(\bx)$ has bounded gradients then it is easy to show that the scaled $\Hdiv$ norm $\|D\cdot\|_{\Hdiv}$ is equivalent to $\|\cdot\|_{\Hdiv}$, \eqref{hdivnorm}. Therefore the bilinear form $\Ac(\cdot, \cdot)$ is continuous and coercive with respect to $\|\cdot\|_{\Hdiv}$ norm, which plays the role of the corresponding energy norm for \eqref{elliptic}. 
\end{remark}
\subsection{Other Boussinesq-Peregrine type systems}
We consider now a modification of Peregrine's system which simplifies the incomplete-elliptic operator $\Ec$.  The new resulting system comprised two decoupled equations for the two components of the depth averaged velocity. This system also generalizes the respective "classical" Boussinesq system of \cite{BCL2015} for general bottoms.

In what follows we consider characteristic quantities for  typical waves in the Boussinesq regime, in particular a typical wave amplitude $a_0$ and length $\lambda_0$ and a typical depth $D_0$. We will denote the linear wave speed by $c_0=\sqrt{gD_0}$. We also consider the dimensionless variables 
\begin{equation}
\tilde{\bf x}=\frac{{\bf x}}{\lambda_0},\qquad  \tilde{t}= \frac{c_0}{\lambda_0}t, \qquad \tilde{\bf u} = \frac{h_0}{a_0c_0}{\bf u}, \qquad \tilde{\eta}=\frac{\eta}{a_0}, \qquad \tilde{D}=\frac{D}{D_0}\ .
\end{equation} 
Then Peregrine's system \eqref{eq:Peregrin} can be written in the nondimensional and scaled form:
\begin{equation}
 \label{eq:prgnd}
\begin{aligned}
& \pd{\tilde{t}}{\tilde{\eta}} + \nabla\cdot\left((\tilde{D}+\varepsilon \tilde{\eta})\tilde{\bu}  \right) = 0\ ,   \\
& \pd{\tilde{t}}{\tilde{\bu}} - \sigma^2 \frac12 \tilde{D} \pd{\tilde{t}}{\nabla\left(\nabla\cdot(\tilde{D}\tilde{\bu}) \right)} + \sigma^2 \frac16 \tilde{D}^2 \pd{\tilde{t}}{\nabla\left(\nabla\cdot\tilde{\bu} \right)} + \varepsilon \left(\tilde{\bu}\cdot\nabla \right)\tilde{\bu} + \nabla \tilde{\eta} = 0\ , 
\end{aligned}
\end{equation}
where $\varepsilon=a_0/D_0$ and $\sigma=D_0/\lambda_0$. In the Boussinesq regime we have $\varepsilon\approx \sigma^2\ll 1$ while the Stokes (or Ursel) number is $S=\varepsilon/\sigma^2=O(1)$.

The classical derivation of Peregrine's system doesn't take into account the irrotationality condition of the flow, \cite{Per67}, which although is not satisfied exactly, it can be expressed \cite{mits2009} as:
\begin{equation}\label{eq:irrotational}
{\tilde{u}}_{\tilde{y}}={\tilde{v}}_{\tilde{x}}+O(\sigma^2)\ .
\end{equation}
Using the former relation we observe that the dispersive terms of Peregrine's system can be simplified to the following:
\begin{equation}
\label{eq:disp1}\
\begin{aligned}
&\nabla\left(\nabla\cdot{\tilde{\bu}}_{\tilde{t}} \right) = \Delta {\tilde{\bu}}_{\tilde{t}}+O(\sigma^2)\ , \\
&\nabla\left(\nabla\cdot(\tilde{D}{\tilde{\bu}}_{\tilde{t}}) \right) = \tilde{D}\Delta {\tilde{\bu}}_{\tilde{t}}+\nabla(\nabla \tilde{D}\cdot {\tilde{\bu}}_{\tilde{t}})+\nabla {\tilde{D}}\nabla\cdot {\tilde{\bu}}_{\tilde{t}}+O(\sigma^2)\ .
\end{aligned}
\end{equation}
From the second equation of \eqref{eq:prgnd} we observe that 
\begin{equation}\label{eq:loword}
\tilde{\bu}_{\tilde{t}}=-\nabla \tilde{\eta}+O(\varepsilon,\sigma^2)\ .
\end{equation}
Using (\ref{eq:loword}) into the second equation of \eqref{eq:disp1} we obtain the relation
\begin{equation}\label{eq:disp1c}
\nabla\left(\nabla\cdot(\tilde{D}{\tilde{\bu}}_{\tilde{t}}) \right) = \tilde{D}\Delta {\tilde{\bu}}_{\tilde{t}}-\nabla(\nabla \tilde{D}\cdot \nabla \tilde{\eta})-\nabla {\tilde{D}}\Delta \tilde{\eta}+O(\varepsilon,\sigma^2)\ .
\end{equation}
Furthermore, using (\ref{eq:irrotational}) we can rewrite the nonlinear term of the second equation of \eqref{eq:prgnd} in the following form
\begin{equation}\label{eq:nonlin}
\left(\tilde{\bu}\cdot\nabla \right)\tilde{\bu}=\frac12\nabla|\tilde{\bu}|^2 + O(\sigma^2)\ .
\end{equation}
Then the  Peregrine system can be written, with the help of (\ref{eq:disp1}), (\ref{eq:disp1c}) and (\ref{eq:nonlin}),  in the form:
\begin{equation}
\label{eq:prgnda}
\begin{aligned}
& \pd{\tilde{t}}{\tilde{\eta}} + \nabla\cdot\left((\tilde{D}+\varepsilon \tilde{\eta})\tilde{\bu}  \right) = 0\ ,    \\
& \pd{\tilde{t}}{\left(\tilde{\bu} - \sigma^2 \frac{1}{3} \tilde{D}^2\Delta \tilde{\bu}\right)} + \frac12\sigma^2\tilde{D}\left( \nabla(\nabla \tilde{D}\cdot \nabla\tilde{\eta})+\nabla\tilde{D}\Delta\tilde{\eta}\right) + \varepsilon \frac12\nabla|\bu|^2 + \nabla \tilde{\eta} = O(\varepsilon^2,\sigma^4)\ .
\end{aligned}
\end{equation}
Discarding the high-order terms, system \eqref{eq:prgnda} is written in dimensional form
\begin{equation}
 \label{eq:prgndd}
\begin{aligned}
& \pd{t}{\eta} + \nabla\cdot\left((D + \eta)\bu \right) = 0\ ,   \\
& \pd{t}{\left(\bu - \frac{1}{3} D^2\Delta \bu\right)} + \frac12 g D\left( \nabla(\nabla D\cdot \nabla\eta)+\nabla D\Delta\eta\right) +  \frac12\nabla|\bu|^2 + g\nabla \eta = 0\ . 
\end{aligned}
\end{equation}
The last system generalises the classical Boussinesq system derived for horizontal bottom $D(\bx)=D_0$ in \cite{BCL2015} to the case of a general bottom. Moreover, the elliptic operators in the momentum equations of (\ref{eq:prgndd}) are decoupled and resemble that of with flat bottom topography \eqref{prgH0}. Similar Boussinesq-type systems have been derived in \cite{mits2009}. It is noted that \eqref{eq:prgnd} cannot be recovered from the analogous regularized Boussinesq systems derived in \cite{mits2009},  even in one-dimensional case \cite{Senth2016} by choosing appropriate coefficients.

Assuming that  $D$ is very smooth, i.e. $\partial^i_x D=\partial^i_yD=O(\varepsilon)$ for $i>0$, we can simplify further  system \eqref{eq:prgnda} to the following {\em simplified Peregrine system} with smooth bottom topography (see also \cite{Chen03}):
\begin{equation}
\label{eq:prgnds}
\begin{aligned}
& \pd{\tilde{t}}{\tilde{\eta}} + \nabla\cdot\left((\tilde{D}+\varepsilon \tilde{\eta})\tilde{\bu}  \right) = 0\ ,    \\
& \pd{\tilde{t}}{\left(\tilde{\bu} - \sigma^2 \frac{1}{3} \tilde{D}^2\Delta \tilde{\bu}\right)} + \varepsilon \frac12\nabla|\tilde{\bu}|^2 + \nabla \tilde{\eta} = O(\varepsilon^2,\sigma^4)\ , 
\end{aligned}
\end{equation}
or written with dimensional variables
\begin{equation}
\label{eq:prgndsd} 
\begin{aligned}
& \pd{t}{\eta} + \nabla\cdot\left((D + \eta)\bu \right) = 0\ ,   \\
& \pd{t}{\left(\bu - \frac{1}{3} D^2\Delta \bu\right)} +  \frac12\nabla|\bu|^2 + g\nabla \eta = 0\ .
\end{aligned}
\end{equation}
We will call system  \eqref{eq:prgndsd}  the {\em simplified} Peregrine system.
Using again \eqref{eq:disp1}(a) we get the system:
\begin{equation}
 \label{eq:prgndm} 
\begin{aligned}
& \pd{\tilde{t}}{\tilde{\eta}} + \nabla\cdot\left((\tilde{D}+\varepsilon \tilde{\eta})\tilde{\bu}  \right) = 0\ .  \\
& \pd{\tilde{t}}{\left(\tilde{\bu} - \sigma^2 \frac{1}{3} \tilde{D}^2\nabla (\nabla\cdot\tilde{\bu})\right)} + \varepsilon \frac12\nabla|\tilde{\bu}|^2 + \nabla \tilde{\eta} = O(\varepsilon^2,\sigma^4)\ , 
\end{aligned}
\end{equation}
or written with dimensional variables
\begin{equation}
 \label{eq:prgndmd} 
\begin{aligned}
& \pd{t}{\eta} + \nabla\cdot\left((D + \eta)\bu \right) = 0\ .  \\  
& \pd{t}{\left(\bu - \frac{1}{3} D^2\nabla(\nabla\cdot \bu)\right)} +  \frac12\nabla|\bu|^2 + g\nabla \eta = 0\ , 
\end{aligned}
\end{equation}
The system in \eqref{eq:prgndmd} will be called {\em modified} Peregrine system. 
It is noted that the differences in the dispersive terms between the simplified and the modified systems required different treatment in the boundary conditions as well. The modified and classical Peregrine systems share the same type of elliptic operator while they have different nonlinear terms. 
\begin{remark}
Existence and uniqueness of solutions to \eqref{eq:Peregrin} is studied in \cite{DuchIsr2018} in the case where the domain is  the whole space. In particular it's been shown that under high-regularity assumptions on the bottom topography $D\in \dot{H}^{N+2}(\mathbb{R}^2)$ for $N\geq 4$, system \eqref{eq:Peregrin} has a local in time, unique solution $(\bu,\eta)\in C([0,T], X^N(\mathbb{R}^2)\times H^N(\mathbb{R}^2))$, where $\dot{H}^N(\mathbb{R}^2)=\{f\in L_{\rm{loc}}^2(\mathbb{R}^2), \nabla f\in H^{N-1}(\mathbb{R}^2)\}$ and $X^N(\mathbb{R}^2)=\{\bu\in L^2(\mathbb{R}^2), \|\bu\|_{X^N}^2=\sum_{|\alpha|=0}^N\|\partial^\alpha \bu\|^2_{L^2}+\mu\|\partial^\alpha \nabla\cdot \bu\|_{L^2}^2<\infty\}$ for $0<\mu\ll 1$. Well-posedness of system \eqref{eq:Peregrin} in bounded domains has only been established  in one-space dimension, cf. \cite{Schonbek1981, Adamy2011}.
\end{remark}

\section{A new BBM-BBM system}

In Peregrine-type systems the mass conservation equation is represented by a first-order partial differential equation. It is known from the theory of finite element method that these kind of differential equations can result to suboptimal convergence rates.  The same behaviour is expected in the case of Peregrine's system as well, \cite{AD2013}. For this reason we proceed also to the derivation of a BBM-type system with small bottom variations, where the mass conservation will contain a regularisation term. 

In \cite{mits2009} a BBM-type system describing weakly nonlinear and weakly dispersive waves propagating over a variable bottom was derived by modifying  both the mass and momentum conservation equations so as to avoid the ${\rm grad-div}$ operator. Here we will perform the same technique but only for the mass equation. We first evaluate the horizontal velocity of the fluid ${\bf u}^\theta$ at some height $z=-D+\theta(\varepsilon\eta+D)$, with $\theta\in[0,1]$. 
Then the velocity ${\bf u}^\theta$ is related to the depth average velocity ${\bf u}$ \cite{Nwogu93} with the formula
\begin{equation}\label{eq:veltheta}
{\bf u}^\theta={\bf u}-\sigma^2\tilde{a}\frac{D}{2}\nabla(\nabla\cdot (D{\bf u}))-\sigma^2\tilde{b}\frac{D^2}{3} \nabla(\nabla\cdot {\bf u})+O(\sigma^4)\ ,
\end{equation}  
with $\tilde{a}=\theta-1/2$ and $\tilde{b}=1/2[(\theta-1)^2-1/3]$. Substitution of (\ref{eq:veltheta}) into \eqref{eq:prgnd} leads to Nwogu's  system
\begin{equation}
 \label{eq:nper1}
\begin{aligned}
& \eta_t+\nabla\cdot((D+\varepsilon\eta){\bf u}^\theta)+\sigma^2\nabla\cdot [\tilde{a}D^2 \nabla(\nabla\cdot(D{\bf u}^\theta))+\tilde{b}D^3\nabla(\nabla\cdot{\bf u}^\theta)]=O(\varepsilon\sigma^2,\sigma^4)\ ,\\
& {\bf u}^\theta_t+\nabla\eta+\varepsilon ({\bf u}^\theta\cdot\nabla){\bf u}^\theta+\sigma^2[\tilde{c}D\nabla(\nabla\cdot(D{\bf u}^\theta_t))+\tilde{d}D^2\nabla(\nabla\cdot{\bf u}^\theta_t)]=O(\varepsilon\sigma^2,\sigma^4)\ , 
\end{aligned}
\end{equation}
with 
\begin{equation}\label{eq:coefsbous}
\tilde{a}=\theta-1/2,~ \tilde{b}=1/2[(\theta-1)^2-1/3],~ \tilde{c}=\theta-1~ \mbox{ and }~\tilde{d}=1/2(\theta-1)^2\ ,
\end{equation}
(see also \cite{mits2009,Nwogu93}). Following again \cite{mits2009}, we have that $D\nabla (\nabla\cdot {\bf u})^\theta=\nabla(\nabla \cdot(D{\bf u}^\theta))-\nabla(\nabla D\cdot {\bf u}^\theta)-\nabla D\nabla\cdot {\bf u}^\theta$  and substitution into \eqref{eq:nper1} we obtain the regularized system
\begin{equation}
\label{eq:nper02}
\begin{aligned}
& \eta_t+\nabla\cdot((D+\varepsilon\eta){\bf u}^\theta)-\sigma^2\tilde{b}\nabla\cdot[D^2(\nabla(\nabla D\cdot {\bf u}^\theta)+\nabla D\nabla\cdot {\bf u}^\theta)]+ \sigma^2(\tilde{a}+\tilde{b})\nabla\cdot [D^2 \nabla(\nabla\cdot(D{\bf u}^\theta))]=O(\varepsilon\sigma^2,\sigma^4)\ , \\
& {\bf u}^\theta_t+\nabla\eta+\varepsilon ({\bf u}^\theta\cdot\nabla){\bf u}^\theta+\sigma^2[\tilde{c}D\nabla(\nabla\cdot(D{\bf u}^\theta_t))+\tilde{d}D^2\nabla(\nabla\cdot{\bf u}^\theta_t)]=O(\varepsilon\sigma^2,\sigma^4)\ .
\end{aligned}
\end{equation}
Observing now from \eqref{eq:nper02} that $\nabla\cdot (D{\bf u}^\theta)=-\eta_t+O(\varepsilon,\sigma^2)$ we write Nwogu's system in a new regularized form as
\begin{equation}
\label{eq:nper03}
\begin{aligned}
& \eta_t+\nabla\cdot((D+\varepsilon\eta){\bf u}^\theta)-\sigma^2\tilde{b}\nabla\cdot[D^2(\nabla(\nabla D\cdot {\bf u}^\theta)+\nabla D\nabla\cdot {\bf u}^\theta)]- \sigma^2(\tilde{a}+\tilde{b})\nabla\cdot [D^2 \nabla \eta_t]=O(\varepsilon\sigma^2,\sigma^4)\ , \\
& {\bf u}^\theta_t+\nabla\eta+\varepsilon ({\bf u}^\theta\cdot\nabla){\bf u}^\theta+\sigma^2[\tilde{c}D\nabla(\nabla\cdot(D{\bf u}^\theta_t))+\tilde{d}D^2\nabla(\nabla\cdot{\bf u}^\theta_t)]=O(\varepsilon\sigma^2,\sigma^4)\ .
\end{aligned}
\end{equation}
To simplify further the system we make the assumption of small bottom variations. For example we assume that the bottom is very smooth in the sense that the derivatives of $D$ are of $O(\varepsilon)$. Although this assumption can result in a system without any contribution of the bottom variations to the dispersive terms, \cite{Chen03}, we choose to simplify only the terms  $\nabla(\nabla D\cdot {\bf u}^\theta)+\nabla D\nabla\cdot {\bf u}^\theta=O(\varepsilon)$ so as to keep as many bottom effects as possible in the dispersive terms and at the same time derive a simple and easy to handle system. The system \eqref{eq:nper03} is simplified then to 
\begin{equation}
\label{eq:nper3}
\begin{aligned}
& \eta_t+\nabla\cdot((D+\varepsilon\eta){\bf u}^\theta)-\sigma^2(\tilde{a}+\tilde{b})\nabla\cdot [D^2 \nabla\eta_t]=O(\varepsilon\sigma^2,\sigma^4)\ , \\
& {\bf u}^\theta_t+\nabla\eta+\varepsilon ({\bf u}^\theta\cdot\nabla){\bf u}^\theta+\sigma^2[\tilde{c}D\nabla(\nabla\cdot(D{\bf u}^\theta_t))+\tilde{d}D^2\nabla(\nabla\cdot{\bf u}^\theta_t)]=O(\varepsilon\sigma^2,\sigma^4)\ .
\end{aligned}
\end{equation}
This is a BBM-BBM-type system because of the existence of the regularized term in the mass conservation equation. A similar system with additional simplifications to the momentum equations was derived in \cite{mits2009}. It is noted that although we can simplify the dispersive terms in the momentum equation as well, we keep the original form because it is easier to use slip boundary conditions (for example this is harder to do with the relatively more simplified system of \cite{mits2009}) and also to be in agreement with Peregrine's and Nwogu's systems. After dropping the high-order terms, we write the system \eqref{eq:nper3} in dimensional form as
\begin{equation}
\label{eq:nper4}
\begin{aligned}
& \eta_t+\nabla\cdot((D+\eta){\bf u}^\theta)-(\tilde{a}+\tilde{b})\nabla\cdot [D^2 \nabla\eta_t]=0\ , \\
& {\bf u}^\theta_t+g\nabla\eta+ ({\bf u}^\theta\cdot\nabla){\bf u}^\theta+\tilde{c}D\nabla(\nabla\cdot(D{\bf u}^\theta_t))+\tilde{d}D^2\nabla(\nabla\cdot{\bf u}^\theta_t)=0\ .
\end{aligned}
\end{equation}
In order for the system \eqref{eq:nper4} to be well posed it is required that $1/3\leq \theta^2\leq 1$, \cite{BCS2002,BCL2015}. Comparing the linear dispersion relation of the system \eqref{eq:nper4} with the corresponding relation of the Euler equations we observe that the parameter $\theta=1$ results to the system which has the closest linear relationship with that of the Euler equations. In order to obtain an exact generalisation of the BBM-BBM system of \cite{BCL2015,BCS2002} we take $\theta=\sqrt{2/3}$ so that $\tilde{a}+\tilde{b}=-(\tilde{c}+\tilde{d})=1/6$ and the system after dropping the $\theta$ in the notation is written 
\begin{equation}
\label{eq:bbm}
\begin{aligned}
& \eta_t+\nabla\cdot((D+\eta){\bf u})-\frac{1}{6}\nabla\cdot \left(D^2 \nabla\eta_t\right)=0\ , \\
& {\bf u}_t+g\nabla\eta+ ({\bf u}\cdot\nabla){\bf u}+\left(\sqrt{\frac{2}{3}}-1\right)D\nabla(\nabla\cdot(D{\bf u}_t))+\left(\frac{5}{6}-\sqrt{\frac{2}{3}}\right)D^2\nabla(\nabla\cdot{\bf u}_t)=0\ .
\end{aligned}
\end{equation}
For the rest of the paper we will call \eqref{eq:bbm}  BBM-BBM system since it coincides with the well-known one-dimension BBM-BBM system when assuming the bottom is flat \cite{BCS2002}. If we take $\theta=\sqrt{1/3}$ then we recover the classical Boussinesq system of \cite{BCS2002} but with different coefficients in front of the dispersive terms compared to those of the classical Peregrine system when the bottom is not flat. Because the system with $\theta=\sqrt{1/3}$ doesn't contribute to the current work,  other than the information that in the Boussinesq regime is reasonable to assume that $\nabla D=O(\varepsilon)$, we will not consider it further.

In order to simulate slip wall boundary conditions with the BBM-BBM system we impose zero Neumann boundary conditions $\nabla\eta\cdot\bn =0$ on the free-surface elevation and $\bu\cdot\bn=0$ on the velocity vector field at the boundary of the computational domain. The boundary condition for $\eta$ is not restrictive but is consistent with the total reflection of water waves on the vertical wall,  \cite{MSM2017}, and it is also satisfied by the Euler equations exactly \cite{Khakimzyanov2018a}. It also appears naturally into the numerical scheme, since it doesn't require any modification of the finite element spaces (and the Sobolev spaces) from the respective spaces we use for the Peregrine system.

\begin{figure}[ht!]
  \centering
  \includegraphics[width=\columnwidth]{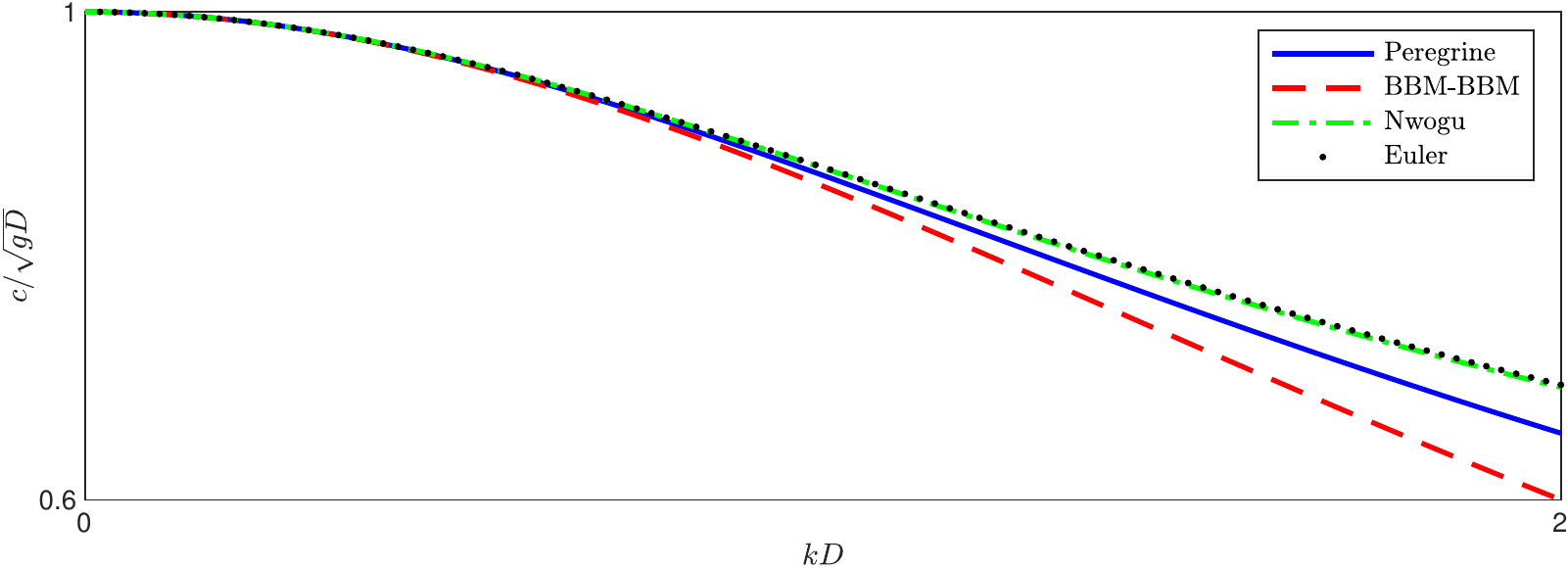}
  \caption{Comparison of linear dispersion relations}
  \label{fig:disprel}
\end{figure}

The linear dispersion relations of the previously mentioned systems are given by the relations 
$$\frac{c_{BBM}}{\sqrt{gD}}=\frac{1}{1+(Dk)^2/6}\quad \mbox{ and }\quad
\frac{c_{Peregrine}}{\sqrt{gD}}=\frac{1}{\sqrt{1+(Dk)^2/3}}\ ,$$
while the linear dispersion relation of the Euler equations is
$$\frac{c_{Euler}}{\sqrt{gD}}=\sqrt{\frac{{\rm tanh}(Dk)}{Dk}}\ .$$
In the previous relations $D$ is assumed to be the depth of the flat bottom and $k$ is the wavenumber of the linear waves.
A comparison of the two dispersion relations of the linearized BBM-BBM and Peregrine systems with those of the linearized Euler equations and Nwogu's system are presented in Figure \ref{fig:disprel}.
It is easily observed that Nwogu's system has better dispersive properties compared to other Boussinesq systems (as expected) while Peregrine system is better than the BBM-BBM system on this account. Of course in the long-wave limit all models converge to the correct value and thus we do not expect discrepancies between the solution of the BBM-BBM system and the other Boussinesq systems. 

\begin{remark}
Following \cite{mits2009} one can derive also the previously mentioned Boussinesq systems \eqref{eq:nper4} assuming that the bottom depends on time, i.e. it is of the form $D(x,y)+\zeta(x,y,t)$, where the time-dependent part is of the same order of magnitude as the free-surface elevation $\eta$. For example the new BBM-BBM system takes the form:
\begin{equation}
\label{eq:bbmvb}
\begin{aligned}
& \eta_t+\nabla\cdot((D+\zeta+\eta){\bf u}^\theta)-(\tilde{a}+\tilde{b})\nabla\cdot [D^2 \nabla\eta_t]=-\tilde{a}\nabla\cdot(D^2\nabla\zeta_t)-\zeta_t\ , \\
& {\bf u}^\theta_t+g\nabla\eta+ ({\bf u}^\theta\cdot\nabla){\bf u}^\theta+\tilde{c}D\nabla(\nabla\cdot(D{\bf u}^\theta_t))+\tilde{d}D^2\nabla(\nabla\cdot{\bf u}^\theta_t)=-\tilde{c}D\nabla\zeta_{tt},
\end{aligned}
\end{equation}
where $\tilde{a}$, $\tilde{b}$, $\tilde{c}$ and $\tilde{d}$ are given by (\ref{eq:coefsbous}). The corresponding Peregrine system (\ref{eq:Peregrin}) can be written as
\begin{equation}
\label{eq:Peregrinvb}
\begin{aligned}
&\eta_t+\nabla\cdot[(D+\zeta+\eta)\bu]=-\zeta_t\ , \\
&{\bf u}_t+g\nabla\eta+(\bu\cdot \nabla\cdot)\bu -\frac{1}{2}D\nabla(\nabla\cdot(D\bu_t))+\frac{1}{6}D^2\nabla(\nabla\cdot\bu_t)=\frac{D}{2}\nabla\zeta_{tt}\ .
\end{aligned}
\end{equation}
The forcing terms of the moving bottom terms can be used not only for simulating waves generated by moving bottoms such as tsunamis \cite{mits2009} but also for the accurate simulation of waves generated due to a wavemaker. We discuss this application in detail in a later section (see also \cite{WB1999,WB2002}).
\end{remark}

All the aforementioned Boussinesq-Peregrine systems  are augmented with initial and boundary conditions and 
can be written in the following compact form, where the values of $\Ec_s, \ \Ec$ and $\Nc$ are presented in Table \ref{PBmodels}
\begin{equation}
\label{PBsys}
\begin{aligned}
 \pd{t}{\Ec_s(\eta)} + \nabla\cdot\left((D+\eta)\bu  \right) & = 0, \\
 \pd{t}{\Ec(\bu)} + \Nc(\bu) + g\nabla\eta & = 0, \\
 \bu \cdot \bn &= 0, \ \text{on} \ \dOmega \ \text{ Original, Modified and BBM-BBM},  \\ 
 \bu & = 0, \  \text{on} \ \dOmega \ \text{ Simplified},  \\
 \nabla\eta\cdot\bn &= 0,  \ \text{on} \ \dOmega \  \text{ BBM-BBM},  \\
 \bu(\bx, 0)  = \bu_0(\bx), \ \eta(\bx,0) &=\eta_0(\bx) .
\end{aligned}
\end{equation} 
 \begin{table}[ht]
 \renewcommand{\arraystretch}{1.6}
\caption{Peregrine-Boussinesq type models}
\centering
\small\addtolength{\tabcolsep}{-3pt}
\begin{tabular}{||c||c||c||c||} \hline
 Model & $\Ec_s(\eta)$ & $\Ec(\bu)$ & $\Nc(\bu)$ \\ \hline\hline
Classical & $\eta$ & $\bu - \frac12 D \nabla\left(\nabla\cdot(D\bu) \right) + \frac16 D^2 \nabla\left(\nabla\cdot\bu \right)$  & $\left(\bu\cdot\nabla \right)\bu$ \\  \hline
Simplified & $\eta$ & $\bu - \frac13 D^2 \Delta\bu$ & $\frac12\nabla|\bu|^2$ \\  \hline
Modified & $\eta$ & $\bu - \frac13 D^2 \nabla(\nabla\cdot\bu)$ & $\frac12\nabla|\bu|^2$ \\  \hline
BBM-BBM & $\eta - \frac{1}{6}\nabla\cdot \left(D^2 \nabla\eta\right)$ & $\bu + \left(\sqrt{\frac{2}{3}}-1\right)D\nabla(\nabla\cdot(D\bu))+\left(\frac{5}{6}-\sqrt{\frac{2}{3}}\right)D^2\nabla(\nabla\cdot\bu)$  & $\left(\bu\cdot\nabla \right)\bu$ \\  \hline
\end{tabular}
\label{PBmodels}
\end{table}%
%
%
%
%
\section{Finite element approximations}\label{sec:feapprox}
We consider now the application of the finite element method for obtaining approximations to the models presented in the previous section. We start with the elliptic problem \eqref{elliptic} and we proceed with the classical Peregrine system \eqref{eq:Peregrin} and its simplifications \eqref{eq:prgndsd}  and \eqref{eq:prgndmd}.

To simplify the presentation of the finite element method we assume in the sequel, without loss of generality, that  $\Omega$ is a polygonal domain that might contain holes and with boundary $\Gamma$.  On $\Omega$ we consider a  regular triangulation $\Th{h}$, \cite{Ciarlet77} consisting of  triangles $K, \ \bar{\Omega}=\cup _K K$, with $h_K$ denoting the diameter of $K$ and $h$ the maximum diameter of these triangles, $h=\max\{h_K, K\in\Th{h}\}$. Further, for a positive integer $r$ let $P^r(K)$ denotes the space of bivariate polynomials of degree $\le r$ on $K$. 

For the approximation of $(\eta, \bu)$ we consider the standard finite element spaces $V_h^r$ consisting of Lagrange-type polynomial basis functions of degree at most $r$:
\begin{equation}
\label{femspace}
V_h^r =\{ \phi \in C(\Omega) : \phi |_K \in P^r(K), K\in\Th{h} \}, \quad U_h^r = V_h^r \times V_h^r , 
\end{equation} 
and $V_h^r$ and $U_h^r$  are finite dimensional subspaces of $V$ and $V\times V$ respectively. 
\begin{remark}
The form of the elliptic part of the Peregrine system \eqref{eq:Peregrin} suggests of using finite element subspaces of  $H(\text{div}, \Omega)$. This approach can be
worked out possibly in the framework of a mixed finite element method where an appropriate discretization of the nonlinear term in \eqref{eq:Peregrin} must be considered. Here, we follow a different approach and we  consider finite element subspaces of $H^1(\Omega)$. 
\end{remark}
\begin{remark} In general we will seek approximations of $\eta$ and $\bu$ from different finite element subspaces $V_h^{r_1}$ and $U_h^{r_2}$  respectively, where $r_1, r_2$ are, possibly different, positive integers. 
\end{remark}
\begin{remark}
The numerical results reported in this and the following sections were obtained by using in part the \emph{FEniCS} computational framework, \cite{Fenics, LMW}. 
\end{remark}

\subsection{The elliptic problem}
We consider first the finite element discretization of the elliptic problem \eqref{elliptic}. An immediate difficult arises due to the boundary condition \eqref{elliptic} which is not satisfied by the elements of $U_h^r$.  To enforce such a  condition on the finite element space $U_h^r$ we follow the well-known technique  of Nitsche (also known to as the Nitsche method), \cite{Nitsche} by introducing an appropriate penalty term in the variational formulation of the problem : we seek $\bu_h \in U_h^r$ such that 
\begin{equation}
\label{femelliptic}
\begin{aligned}
\Cc(\bu_h, \bphi):= \Ac(\bu_h, \bphi) + & \Bc(\bu_h,\bphi) + \frac13 \frac{C_N}{h}\langle D\bu_h\cdot \bn, D\bphi\cdot\bn \rangle = (\bF, \bphi), \quad \forall \bphi\in U_h^r  , \\
\text{ where } \quad \Bc(\bu_h, \bphi) = & -\frac13 \intg \Div(D\bu_h)\, D\bphi\cdot\bn\, \ds - \frac16 \intg  (\bu_h\cdot\nabla D)\, D\bphi\cdot\bn\, \ds \\
& -\frac13 \intg \Div(D\bphi) \, D\bu_h\cdot\bn\, \ds - \frac16 \intg (\bphi\cdot\nabla D)\, D\bu_h\cdot\bn\, \ds , 
\end{aligned} 
\end{equation}
where $\bn$ is the outward unit normal to the boundary $\Gamma$ and  $C_N$ is constant to be chosen.  The bilinear form $\Ac(\cdot, \cdot)$ is defined as in \eqref{weakelliptic1}. The first two terms in $\Bc(\cdot, \cdot)$ are from integration by parts while the last two terms are added for symmetry and they are consistent in the sense that they vanish in the continuous case due to the boundary condition of \eqref{elliptic}. Along the lines of \eqref{weakelliptic2} one can show that $\Cc(\cdot, \cdot)$ is continuous and coercive with respect to the $\|D\cdot\|_{\Hdiv}$, so \eqref{femelliptic} has unique solution. To validate the method we verify the experimental order of convergence(EOC) by means of the method of
manufactured solutions. To this end, from \eqref{elliptic}  we compute an   $\bF$ such that 
\begin{equation}
\label{exactu}
\bu(\bx) = (\cos(\frac{\pi y}{2}) \sin(\pi x), \cos(\frac{\pi x}{2}) \sin(\pi y)), \quad (x,y)\in [0,1]\times[0,1] , 
\end{equation}
is an exact solution that  satisfies the boundary condition $\bu\cdot\bn = 0$.  We compute the $EOC$'s based on the error between the exact and approximate solution measured in the $L^2, \ H^1$ and $H(\text{div})$ norms:
\begin{equation}
\label{EOC}
EOC = \log\left(\frac{E_1}{E_2}\right)/ \log\left(\frac{h_1}{h_2}\right), 
\end{equation}
where $E_1, E_2$ are the normed errors and $h_1, h_2$ are the corresponding mesh sizes. First we consider the case with flat bottom $D(\bx)\equiv 1$. 
In Tables \ref{uP1}, \ref{uP2} the convergence rates of $\bu$ are presented for various values of $r$ and $C_N=50$. There were no essential difference observed in the errors or rates by taking larger values of penalty constant $C_N$. The optimal convergence rates are observed in all norms. 
\begin{table}[ht]
  \centering
  \caption{Convergence rates for $\bu$: flat bottom,  $r=1, \ C_N=50$}
    \begin{tabular}{||c||c|c||c|c||c|c||}\hline
    h & $\|\bu-\bu_h\|$ &  \text{EOC} &  $\|\bu-\bu_h\|_1$ &  \text{EOC} & $\|\bu-\bu_h\|_{\Hdiv}$ & \text{EOC} \\ \hline
  1.250e-01 &  4.564e-03 &   &  2.481e-01 &   &  1.845e-01 &  \\ \hline
  8.333e-02 &  2.027e-03 &  2.001 &  1.655e-01 &  0.998 &  1.230e-01 &  1.000 \\ \hline
  6.250e-02 &  1.140e-03 &  2.001 &  1.241e-01 &  0.999 &  9.225e-02 &  1.000 \\ \hline
  5.000e-02 &  7.297e-04 &  2.000 &  9.933e-02 &  1.000 &  7.380e-02 &  1.000 \\ \hline
  4.167e-02 &  5.067e-04 &  2.000 &  8.278e-02 &  1.000 &  6.150e-02 &  1.000 \\ \hline
  3.571e-02 &  3.722e-04 &  2.000 &  7.095e-02 &  1.000 &  5.271e-02 &  1.000 \\ \hline
  3.125e-02 &  2.850e-04 &  2.000 &  6.208e-02 &  1.000 &  4.612e-02 &  1.000 \\ \hline
\end{tabular}
  \label{uP1}%
\end{table}%
\begin{table}[htbp]
  \centering
  \caption{Convergence rates for $\bu$: flat bottom, $r=2, \ C_N=50$}
    \begin{tabular}{||c||c|c||c|c||c|c||}\hline
    h & $\|\bu-\bu_h\|$ &  \text{EOC} &  $\|\bu-\bu_h\|_1$ &  \text{EOC} & $\|\bu-\bu_h\|_{\Hdiv}$ & \text{EOC} \\ \hline  
  1.250e-01 &  1.290e-04 &   &  1.059e-02 &   &  7.185e-03 &   \\ \hline  
  8.333e-02 &  3.901e-05 &  2.950 &  4.725e-03 &  1.990 &  3.194e-03 &  2.000 \\ \hline  
  6.250e-02 &  1.661e-05 &  2.968 &  2.662e-03 &  1.994 &  1.796e-03 &  2.000 \\ \hline  
  5.000e-02 &  8.548e-06 &  2.977 &  1.706e-03 &  1.996 &  1.150e-03 &  2.000 \\ \hline  
  4.167e-02 &  4.963e-06 &  2.982 &  1.185e-03 &  1.997 &  7.984e-04 &  2.000 \\ \hline  
  3.571e-02 &  3.133e-06 &  2.985 &  8.711e-04 &  1.997 &  5.865e-04 &  2.000 \\ \hline  
  3.125e-02 &  2.102e-06 &  2.987 &  6.671e-04 &  1.998 &  4.490e-04 &  2.000 \\ \hline  
\end{tabular}
  \label{uP2}%
\end{table}%
The case of a non-flat bottom is different. We take $D(\bx) = -\frac{1}{20}(x+y) + \frac32$ and we compute the EOC in this case. The results are presented in Tables \ref{uP1nf} and \ref{uP2nf}. We note that for linear finite elements, $r=1$, the optimal convergence rates are observed in all norms. However, for $r=2$ the convergence rates are suboptimal for the $L^2$ and $H^1$ norms but optimal in the $H(\text{div})$ norm. The same behaviour is observed for $r=3$. 

\begin{table}[htbp]
  \centering
  \caption{Convergence rates for $\bu$: non-flat bottom,  $r=1, \ C_N=50$}
    \begin{tabular}{||c||c|c||c|c||c|c||}\hline
    h & $\|\bu-\bu_h\|$ &  \text{EOC} &  $\|\bu-\bu_h\|_1$ &  \text{EOC} & $\|\bu-\bu_h\|_{\Hdiv}$ & \text{EOC} \\ \hline
  1.250e-01 &  4.789e-03 &   &  2.480e-01 &   &  1.845e-01 &   \\ \hline
  8.333e-02 &  2.128e-03 &  2.000 &  1.655e-01 &  0.998 &  1.230e-01 &  1.000 \\ \hline
  6.250e-02 &  1.197e-03 &  2.000 &  1.242e-01 &  0.999 &  9.225e-02 &  1.000 \\ \hline
  5.000e-02 &  7.664e-04 &  1.999 &  9.933e-02 &  0.999 &  7.380e-02 &  1.000 \\ \hline
  4.167e-02 &  5.323e-04 &  1.999 &  8.278e-02 &  1.000 &  6.150e-02 &  1.000 \\ \hline
  3.571e-02 &  3.911e-04 &  1.999 &  7.096e-02 &  1.000 &  5.271e-02 &  1.000 \\ \hline
  3.125e-02 &  2.995e-04 &  1.999 &  6.209e-02 &  1.000 &  4.612e-02 &  1.000 \\ \hline
\end{tabular}
  \label{uP1nf}%
\end{table}%
\begin{table}[htbp]
  \centering
  \caption{Convergence rates for $\bu$: non-flat bottom, $r=2, \ C_N=50$}
    \begin{tabular}{||c||c|c||c|c||c|c||}\hline
    h & $\|\bu-\bu_h\|$ &  \text{EOC} &  $\|\bu-\bu_h\|_1$ &  \text{EOC} & $\|\bu-\bu_h\|_{\Hdiv}$ & \text{EOC} \\ \hline  
  1.250e-01 &  1.425e-04 &   &  1.186e-02 &   &  7.177e-03 &   \\ \hline
  8.333e-02 &  4.799e-05 &  2.684 &  5.922e-03 &  1.714 &  3.191e-03 &  1.999 \\ \hline
  6.250e-02 &  2.311e-05 &  2.540 &  3.788e-03 &  1.553 &  1.796e-03 &  1.999 \\ \hline
  5.000e-02 &  1.348e-05 &  2.416 &  2.758e-03 &  1.421 &  1.149e-03 &  2.000 \\ \hline
  4.167e-02 &  8.831e-06 &  2.321 &  2.168e-03 &  1.321 &  7.981e-04 &  2.000 \\ \hline
  3.571e-02 &  6.243e-06 &  2.250 &  1.788e-03 &  1.248 &  5.864e-04 &  2.000 \\ \hline
  3.125e-02 &  4.655e-06 &  2.198 &  1.525e-03 &  1.195 &  4.489e-04 &  2.000 \\ \hline
\end{tabular}
  \label{uP2nf}%
\end{table}%
\subsection{Fully discrete schemes}
We consider now the finite element discretization of systems \eqref{PBsys}.  For $t>0$,  we seek $(\eta_h(t), \bu_h(t)) \in V_h^{r_1}\times U_h^{r_2}$ such that 
\begin{equation}
\begin{aligned}
(\pd{t}{\Ec_s(\eta_h)},\psi) - F_s(\eta_h, \bu_h, \psi)  &=  0, \quad  F_s(\eta_h, \bu_h, \psi) = - (\nabla\cdot\left((D+\eta_h)\bu_h\right),\psi), \quad \forall \psi\in V_h^{r_1} ,  \\
(\pd{t}{\Ec(\bu_h)}, \bphi)  - F_u(\eta_h, \bu_h, \bphi) & = 0, \quad  F_u(\eta_h, \bu_h, \bphi) = -  (\Nc(\bu_h),\bphi) - g(\nabla\eta_h, \bphi), \quad \forall \bphi\in U_h^{r_2} . 
\end{aligned}
\end{equation}
Performing a standard integration by parts, as in the elliptic case \eqref{femelliptic} we obtain 
\begin{equation}
\label{PBfem}
\begin{aligned}
\pd{t}{\Cc_s(\eta_h},\psi)  & = F_s(\eta_h, \bu_h, \psi),  \quad \forall \psi\in V_h^{r_1} , \\
\pd{t}{\Cc(\bu_h}, \bphi) & = F_u(\eta_h, \bu_h, \bphi), \quad \forall \bphi\in U_h^{r_2} , \\
\eta_h(0) = \Pc \eta_0, & \ \bu_h(0) = \Pc \bu_0 ,
\end{aligned}
\end{equation}
where $\Cc(\cdot, \cdot)$ is as in \eqref{femelliptic} and $\Cc_s(\eta_h,\psi)=(\eta_h, \psi)$ or $\Cc_s(\eta_h,\psi)=(\eta_h, \psi) + (\nabla\eta_h, \nabla\psi)$ and $\Pc$ denotes the $L^2$-projection on the finite element space.  The bilinear forms $\Cc(\cdot, \cdot)$ and $\Cc_s(\cdot, \cdot)$ are continuous and coercive with respect to $\|D\cdot\|_{\Hdiv}$ and $\|\cdot\|_1$  respectively, thus  system \eqref{PBfem} has a unique solution. 

The system of ordinary differential equations \eqref{PBfem} is discretized by using the classical fourth-order accurate, four-stage explicit Runge-Kutta method(RK4). Indeed this method is used as a time stepping mechanism for all the numerical experiments presented in the sequel.

\subsection{Spatial convergence rates and stability}
We turn now our attention to the convergence of the fully-discrete scheme.  We begin with the  validation of  the  scheme by performing accuracy tests and computing the EOC's for various choices of the finite element spaces. 

It is known that the standard Galerkin method applied to Peregrine's system in 1D with reflective boundary conditions converges with a suboptimal convergence rate, \cite{AD2012, AD2013}. Specifically, it was shown that for linear elements  the convergence rates of the $L^2$ error for the velocity $\bu$ was optimal $2$, while for the free surface $\eta$ was suboptimal and equal to $3/2$. The convergence rates of the $H^1$ error was also shown to decrease with a suboptimal rate of $1/2$ for the free-surface variable $\eta$, while the convergence rate of the respective $H^1$ error for the velocity variable was optimal and equal to $1$.

In similar situations it has been observed that using different finite element spaces for the dependent variables can be beneficial and can result to a fully-discrete scheme that converges with optimal rate to the exact solution, \cite{AD2010,MSM2017}. On the other hand, in the case of other Boussinesq-type system such as the Bona-Smith system, the order of convergence can be different between one and two space dimensions, \cite{DMS2007,ADM2010}. 
For this reasons, it's unrealistic  to have expectations for specific convergence rates for our numerical method unless we check. 

We  compute the EOC's by the method of manufactured solutions. We choose an exact solution $(\eta, \bu)$,  we compute appropriate right hand side $\bF$ using the differential equations and then solve a non-homogeneous problem. For the classical, modified and the new Peregrine systems we choose an exact solution that satisfies the boundary condition $\bu\cdot\bn=0$. However for  the simplified Boussinesq system we need to choose a different exact solution since  an extra boundary condition $\bu=0$ is needed due to the presence of the Laplacian in the differential operator.  As for $\eta$,  Neumann boundary condition is taken for the BBM-BBM system due to the presence of the BBM-term in the continuity equation, while no extra boundary condition is considered for the other models. The value of $\eta$ is the same for all systems. The exact values for  $(\eta, \bu)$ are :
\begin{equation}
\begin{aligned}
& \bu(\bx, t) = e^{(x+y)t}\left(\cos\left(\frac{\pi}{2}y\right)\sin(\pi x),\ \cos\left(\frac{\pi}{2}x\right)\sin(\pi y)\right), \ \text{ Classical, Modified, BBM-BBM },  \\
& \bu(\bx, t) =e^t \left( x \cos\left(\frac{\pi}{2}x\right)\sin(\pi y), \ y \cos\left(\frac{\pi}{2}y\right)\sin(\pi x)\right), \quad\text{ Simplified },  \\
& \eta(\bx,t) = e^t \cos(\pi x)\sin(\pi y) .
\end{aligned}
\end{equation}
We also consider a linearly varying bottom topography 
$$H(x,y)=-\frac{1}{20}(x+y)+\frac{3}{2}\ ,$$
which represents a sloping bottom with slope $1/20$. We integrate numerically the non-homogeneous all Peregrine-Boussinesq type system and we compute the experimental orders of convergence ($EOC$) based on the $L^2, \ H^1$ and $\Hdiv$ errors between the numerical and the exact solutions, as in (\ref{EOC}). 
All the computations where performed in the unit square $\Omega=[0,1]\times[0,1]$ with final time $T=1$. On $\Omega$ we consider a regular uniform triangulation $\Th{h}$ consisting of equilateral triangles and we denote as $h$ the length of a side. 
The temporal time step was taken $\Delta t=5\times 10^{-4}$ to ensure that the contribution of the time-stepping method to the total error is negligible. For approximation we use  finite element spaces $V_h^{r_1}, \ U_h^{r_2}$ with various polynomial degrees $r_1,r_2=1, 2, 3, 4$. 

\begin{figure}[h!]
  \centering
  \includegraphics[width=18cm]{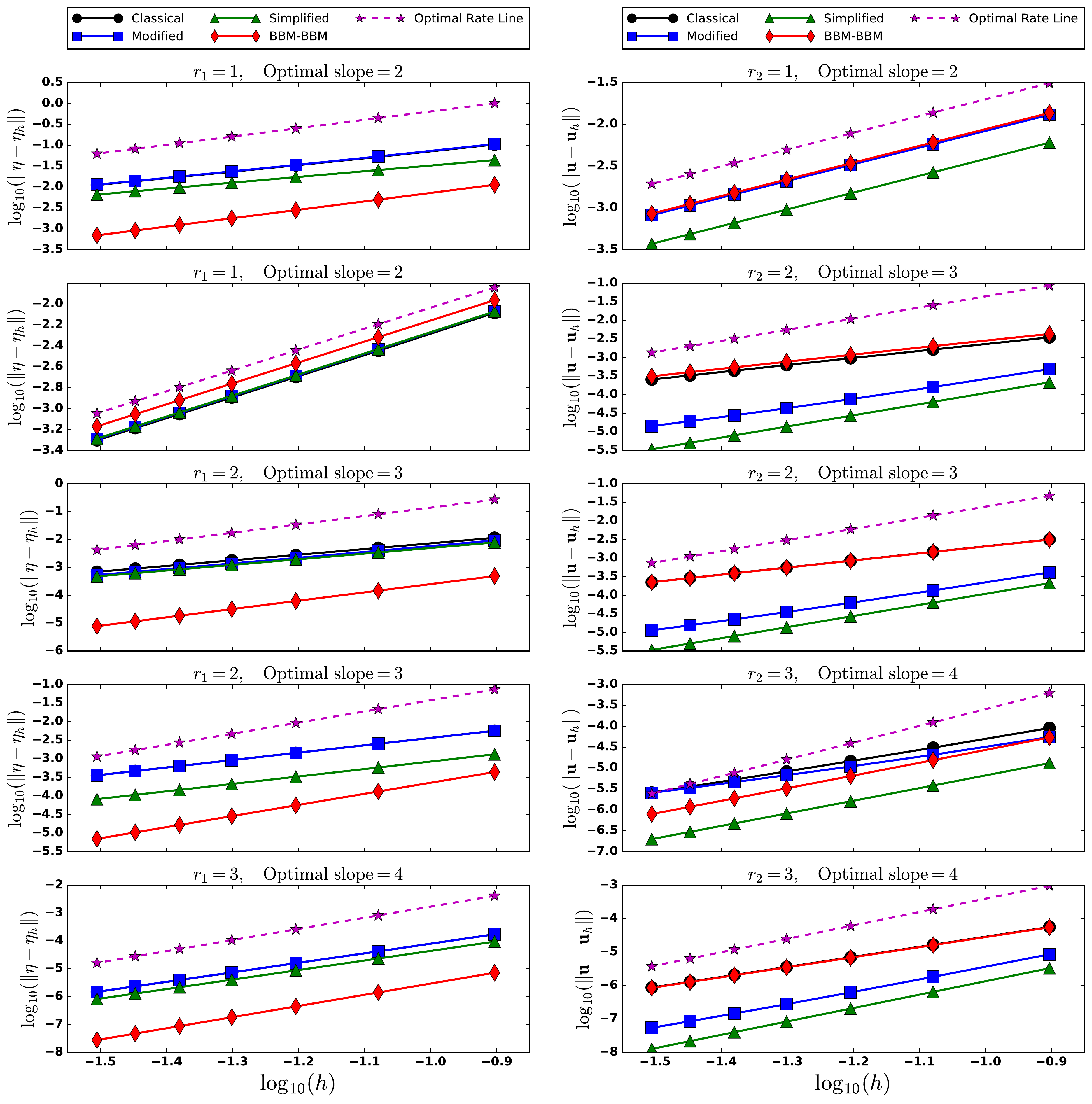}
  \caption{Log-Log plot of errors in the $L^2$-norm : $\eta$ in left column and $\bu$ in right column.}
  \label{Lrates}
\end{figure}

\begin{figure}[h!]
  \centering
  \includegraphics[width=18cm]{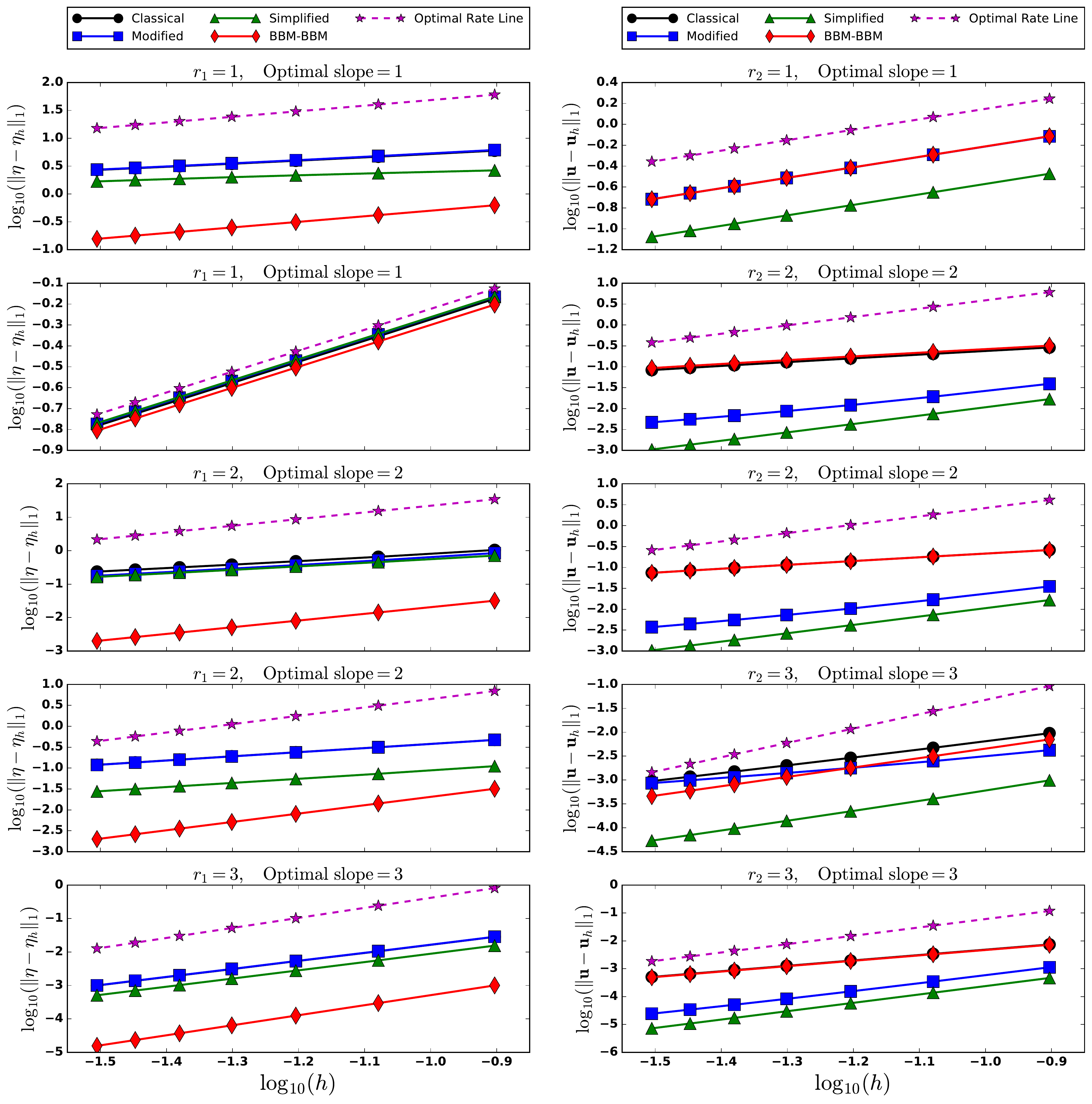}
  \caption{Log-Log plot of errors in the $H^1$-norm : $\eta$ in left column and $\bu$ in right column.}
  \label{Hrates}
\end{figure}

The errors measured in $L^2$ and $H^1$ norms for $(\eta, \bu)$ are shown in Figures \ref{Lrates} and \ref{Hrates} respectively for five different combinations 
of the values $r_1$ and $r_2$.   The behaviour of the error for each variable and each model is presented with solid lines in log-log scale for various values of the discretization parameter $h$. In each graph the corresponding optimal rate is depicted by a dashed line with optimal slope. 

We consider first $(\eta_h, \bu_h)$  in finite element spaces $(V_h^{r}, \ U_h^{r})$ of the same polynomial degree $r$ with $r=1,2,3,4$.  The EOC's for $\eta$ in the $L^2$ and $H^1$ norms for the \emph{classical, modified} and \emph{simplified} systems are exactly the same for all values of $r$, but suboptimal by $1/2$ for the odd degree polynomials and suboptimal by $1$ for the even degree basis functions.  For the \emph{BBM-BMM} system the optimal EOC is obtained in both norms. The latter is attributed to the presence of a BBM-type term in the continuity equation, \eqref{eq:bbm}(a). 
The situation is somehow different for the EOC's for $\bu$.  For $r=1$ the EOC's are optimal in the $L^2$ and $H^1$ norms for all systems. However for $r=2,3,4$ the situation differs and for the \emph{classical, modified} and \emph{BBM-BBM} systems the EOC's in both $L^2$ and $H^1$ norms are suboptimal by $1$. On the other hand, the optimal EOC is achieved for the simplified system since the corresponding elliptic operator involves the Laplacian. 

We tested also three cases where the degree of polynomials $(r_1, r_2)$ of the finite element spaces $(V_h^{r_1}, \ U_h^{r_2})$ they  differ by $1$ and in particular we took $(r_1, r_2) = (1, 2), (2, 3), (3, 4)$. For odd values of $r_1$  the optimal EOC's for $\eta$ were obtained in both norms and for all systems, while for even values of $r_1$ all the EOC's are suboptimal by $1$, except for the BBM-BBM system which the EOC's are optimal.   For $\bu$ for odd values of $r_2$ the resulting EOC's are suboptimal by $1$ in both norms for the \emph{classical, modified} and \emph{BBM-BBM}, but optimal for the \emph{simplified} system. For even values of $r_2$ the convergence rates for $\bu$ are suboptimal by $1$ in both norms for the \emph{simplified, BBM-BBM} systems and suboptimal by $2$ for the \emph{classical} and \emph{modified} systems. Convergence rates for 
$\bu$ were also computed  in the $\Hdiv$ norm. These rates proven to be optimal for all the systems and all the choices of  $r_1$ and $r_2$. 
 
All these results are summarized in Table \ref{EOCs} where by $E_k(\cdot), \ k=0,1$ we denoted the EOC in the $L^2$ and $H^1$ respectively, while $E_d(\cdot)$ is the corresponding EOC in the $\Hdiv$ norm. 

\begin{table}[ht]
\caption{EOC's for $(\eta, \bu)$ for various choices of polynomial degree $r_1, \ r_2$}
\centering
\begin{tabular}{||c||c|c||c|c|c||c|c||c|c|c||}\hline
 & \multicolumn{5}{c||}{Classical} & \multicolumn{5}{c||}{Modified}   \\ \hline
 & \multicolumn{2}{c||}{$\eta$} & \multicolumn{3}{c||}{$\bu$} &  \multicolumn{2}{c||}{$\eta$} & \multicolumn{3}{c||}{$\bu$} \\ \hline
$(r_1,r_2)$ & $E_0(\eta)$ & $E_1(\eta)$ & $E_0(\bu)$  & $E_1(\bu)$  & $E_d(\bu)$ & $E_0(\eta)$ & $E_1(\eta)$ & $E_0(\bu)$  & $E_1(\bu)$  & $E_d(\bu)$\\ \hline
$(1,1)$ & 1.5 & 0.5      & 2 & 1     & 1    & 1.5 & 0.5      & 2 & 1 & 1  \\ \hline
$(2,2)$ & 2 & 1            & 2 & 1     & 2    & 2 & 1            & 2 & 1 & 2  \\ \hline
$(3,3)$ & 3.5 & 2.5      & 3 & 2     & 3    & 3.5 & 2.5      & 3 & 2 & 3  \\ \hline
$(4,4)$ & 4 &  3          & 4 & 3     & 4    & 4 & 3            & 4 & 3 & 4  \\ \hline\hline
$(1,2)$ & 2 & 1            & 2 & 1     & 2    & 2 & 1            & 2 & 1 & 2  \\ \hline
$(2,3)$ & 2.5 & 1         & 2 & 1.5  & 2    & 2 & 1            & 2 & 1 & 3  \\ \hline
$(3,4)$ & 4 & 3            & 4 & 3      & 4   & 4 & 3            & 4 & 3 & 4  \\ \hline\hline        
 & \multicolumn{5}{c||}{Simplified} & \multicolumn{5}{c||}{BBM-BBM}   \\ \hline
 & \multicolumn{2}{c||}{$\eta$} & \multicolumn{3}{c||}{$\bu$} &  \multicolumn{2}{c||}{$\eta$} & \multicolumn{3}{c||}{$\bu$} \\ \hline
$(r_1,r_2)$ & $E_0(\eta)$ & $E_1(\eta)$ & $E_0(\bu)$  & $E_1(\bu)$  & $E_d(\bu)$ & $E_0(\eta)$ & $E_1(\eta)$ & $E_0(\bu)$  & $E_1(\bu)$  & $E_d(\bu)$ \\ \hline
$(1,1)$ & 1.5 & 0.5      & 2 & 1 & 1       & 2 & 1       & 2 & 1 & 1   \\ \hline
$(2,2)$ & 2 & 1            & 3 & 2 & 2       & 3 & 2       & 2 & 1 & 2   \\ \hline
$(3,3)$ & 3.5 & 2.5      & 4 & 3 & 3       & 4 & 3       & 3 & 2 & 3   \\ \hline
$(4,4)$ & 4 & 3            & 5 & 4 & 4       & 5 & 4       & 4 & 3 & 4  \\ \hline\hline
$(1,2)$ & 2 & 1            & 3 & 2 & 2       & 2 & 1       & 2 & 1 & 2  \\ \hline
$(2,3)$ & 2 & 1            & 3 & 2 & 3       & 3 & 2       & 3 & 2 & 3  \\ \hline 
$(3,4)$ & 4 & 3            & 5 & 4 & 4       & 4 & 3       & 4 & 3 & 4  \\ \hline
\end{tabular}
\label{EOCs}
\end{table}%
\begin{remark}
The computational convergence rates are the same as in the one-dimensional case, \cite{MSM2017}, which suggests the same behaviour of the system in two-dimensions.
\end{remark}

The presence of the incomplete-elliptic operator $\Ec(\bu)$, \eqref{elliptic} in the left hand side of the systems \eqref{PBsys} makes the corresponding o.d.e systems \eqref{PBfem} mildly stiff and the fully discrete schemes are stable under a relaxed CFL-type condition. We provide now some indicative numerical results 
concerning the size of the Courant number. To do so, we consider  a solitary wave of various 
amplitudes $A$ and we study its propagation using the classical Peregrine system \eqref{eq:Peregrin} and the newly derived BBM-BBM system \eqref{eq:bbm}. 
For both systems we compute the largest Courant number such that the propagation is stable. We consider a channel $[-50, 100]$ of constant depth $D(\bx)=1$ and the initial location of the solitary wave is at $x=-20$. The channel is covered by an unstructured triangulation of $33192$ triangles with a minimum size $h_{min}\sim0.09$. We integrate both systems up to $T=100$  and we compute Courant number defined as $c_s \frac{\dt}{h_{min}}$, where $c_s= \sqrt{g(D+A)}$ is the phase speed of the solitary wave and  $\dt$ is the time step. We tested solitary waves with five different amplitudes, namely we took $A=0.05, 0.1, 0.2, 0.4, 0.5$. The maximum allowed Courant numbers in all test cases were of the order of $10$ showing that the method is essentially unconditionally stable, \cite{AD2010}.

\section{Solitary waves}\label{sec:solitwaves}

In this section we study the accuracy of the numerical methods in the propagation of a solitary wave and we present some differences in the behaviour of the numerical solutions between the Peregrine and BBM-BBM systems. Since both systems with flat bottom are reduced to the respective Boussinesq systems of \cite{BCS2002} it is known that they possess classical line solitary wave solution \cite{DM2008} that propagate with constant speed and without change in shape at any direction. Due to the lack of exact analytical formulas for the solitary wave solutions of the systems at hand we employed the Petviashvili method for their numerical generation.

\subsection{The Petviashvili iteration}

The Petviashvili method  \cite{Petv1976},  is a modified fixed point method for solving nonlinear equations originally derived from equations for solitary waves. Although Petviashvili method has been applied and studied analytically using Fourier / Pseudospectral methods \cite{PS2004,AD2014} it has been applied successfully with Galerkin finite element methods as well \cite{OSSS2016}.

Here we first present the Petviashvili iteration for the approximation of the line solitary waves of the Peregrine system using a Galerkin finite element method. Assume that we seek for traveling wave solutions of the form 
\begin{equation}\label{eq:ansatz}
\eta(\bx,t)=\eta(\boldsymbol{\xi}),\quad \bu(\bx,t)=\bu(\boldsymbol{\xi})\ ,
\end{equation} where $\boldsymbol{\xi}=(\xi,\zeta)$ with $\xi=\boldsymbol{\alpha} \bx -c_s t-\bx_0$, $\bx_0\in \mathbb{R}$, and $\zeta=\boldsymbol{\alpha}^{\perp}\bx$. This represents a traveling wave solution that travels with constant speed $c_s$, without change in shape,  in a channel aligned along the direction of the vector $\boldsymbol{\alpha}=(\alpha_x,\alpha_y)^T$, $|\boldsymbol{\alpha}|=1$ over a flat bottom $D(\bx)=D_0$, that tends to zero at infinity, while there is no restrictions along the vertical to $\boldsymbol{\alpha}$ direction $\boldsymbol{\alpha}^{\perp}=(\alpha_y,-\alpha_x)^T$. Because the domain for the numerical simulations must be bounded, we assume that the channel is long enough so as the exponentially decreasing solitary wave fits along the $\boldsymbol{\alpha}$ direction in the sense that at the edges of the domain the solitary wave is practically zero.  
Substituting the {\em ansatz} (\ref{eq:ansatz})  into Peregrine's equations (\ref{eq:Peregrin}) we obtain the system of partial differential equations
\begin{equation}\label{eq:sysode1}
\begin{aligned}
&-c_s\eta_\xi+\nabla_{\xi\zeta}\cdot (D_0+\eta) A\bu=0\ ,\\
&-c_s\bu_\xi+gA\nabla_{\xi\zeta}\eta+(\bu\cdot A\nabla_{\xi\zeta})\bu+c_s\frac{1}{3}D_0^2A\nabla_{\xi\zeta}(\nabla_{\xi\zeta}\cdot A\bu_\xi)=0 \ ,
\end{aligned}
\end{equation}
where $\nabla_{\xi\zeta}=(\partial_\xi,\partial_\zeta)^T$ and
$$A=\begin{pmatrix} \alpha_x & \alpha_y\\ \alpha_y & -\alpha_x\end{pmatrix}\ .$$
Multiplying the second equation in (\ref{eq:sysode1}) with $A$ and 
denoting $\bw=(w,\tilde{w})$ where
\begin{equation}\label{eq:petsystem0}
w=\boldsymbol{\alpha}\cdot \bu, \qquad \tilde{w}=\boldsymbol{\alpha}^\perp \cdot \bu\ ,
\end{equation}
then the system (\ref{eq:sysode1}) becomes
\begin{equation}\label{eq:sysode1a}
\begin{aligned}
&-c_s\eta_\xi+\nabla_{\xi\zeta}\cdot (D_0+\eta) \bw=0\ ,\\
&-c_s\bw_\xi+g\nabla_{\xi\zeta}\eta+(\bw\cdot \nabla_{\xi\zeta})\bw+c_s\frac{1}{3}D_0^2\nabla_{\xi\zeta}(\nabla_{\xi\zeta}\cdot \bw_\xi)=0 \ .
\end{aligned}
\end{equation}
Assuming that the solution is constant along the direction of the vector $\boldsymbol{\alpha}^\perp$ then system (\ref{eq:sysode1a}) can be simplified after integration into the system
\begin{equation}\label{eq:sysode1b}
\begin{aligned}
-c_s\eta+ (D_0+\eta) w&=0\ ,\\
-c_sw+g\eta+\frac{1}{2}w^2+c_s\frac{1}{3}D_0^2w_{\xi\xi}&=0 \ , \\
\tilde{w}&=0 \ ,
\end{aligned}
\end{equation}
where we intentionally keep the dependence of the equations on $\tilde{w}$ to emphasize the fact that we are looking for two-dimensional waves and also for implementation purposes.
Eliminating $\eta$ in (\ref{eq:sysode1b}) we obtain the relation 
\begin{equation}\label{eq:reletaw}
\eta=\frac{D_0w}{c_s-w}\ ,
\end{equation}
between the unknowns $\eta$ and $w$ and the partial differential equation
\begin{equation}\label{eq:petv1}
\mathcal{L}\bw=\mathcal{N}(\bw)\ ,
\end{equation}
with 
\begin{equation}
\mathcal{L}=\begin{pmatrix} c_s-\frac{1}{3}D_0^2c_s\partial_{\xi\xi} & 0 \\
0 & 1 \end{pmatrix}~ \mbox{ and } ~ \mathcal{N}(\bw)=\begin{pmatrix} \frac{1}{2}w^2+g\frac{D_0w}{c_s-w} \\ 0\end{pmatrix}\ .
\end{equation}
Solving the equation (\ref{eq:petv1}) for $w$ we recover the free-surface elevation from (\ref{eq:reletaw}). In order to solve (\ref{eq:petv1}) we employ again the standard Galerkin finite element method: seek  $\bw_h\in U_h^{r_2}$ satisfying the equation written in the weak form
\begin{equation}\label{eq:petv2}
\mathcal{L}_h(\bw_h,\boldsymbol{\chi})=(\mathcal{N}(\bw_h),\boldsymbol{\chi}),\quad\mbox{for all}\quad \boldsymbol{\chi}\in U_h^{r_2}, 
\end{equation}
where $\mathcal{L}_h(\bw,\boldsymbol{\chi})$ is the bilinear form 
$$\mathcal{L}_h(\bw,\boldsymbol{\chi})=c_s(w,\phi)+\frac{1}{3}D_0^2c_s(w_\xi,\phi_\xi)+(\tilde{w},\psi), $$ defined for all $\bw=(w,\tilde{w})\in U_h^{r_2}$ and $\boldsymbol{\chi}=(\phi,\psi)\in U_h^{r_2}$. In order to solve the nonlinear equation (\ref{eq:petv2}), we use the Petrviashvilli iteration defined as 
\begin{equation}
\mathcal{L}_h(\bw_h^{n+1},\boldsymbol{\chi})=M_n^\gamma (\mathcal{N}(\bw_h^n),\boldsymbol{\chi}),\quad n=0,1,\cdots\ ,
\end{equation}
for all $\boldsymbol{\chi}\in U_h^{r_2}$ and $M_n$ is the scalar  
\begin{equation}
M_n=\frac{\mathcal{L}_h(\bw^n_h,\bw^n_h)}{(\mathcal{N}(\bw^n_h),\bw^n_h)}\ .
\end{equation}
The function $\bw_h^0$ can be the $L^2$-projection of any initial guess $\bw^0$ of the solitary wave solution $w$ onto the finite element space $V_h^{r_2}$. As an approximation to the line solitary waves of the Boussinesq systems we considered the functions
\begin{equation}\label{eq:initialguesw}
\begin{aligned}
\eta^0(\boldsymbol{\xi})&=A{\rm sech}^2\left(\lambda \xi\right)\ , \\
w^0(\boldsymbol{\xi})&=c_s\eta^0(\boldsymbol{\xi})/(D_0+\eta^0(\boldsymbol{\xi}))\ , \\
\tilde{w}^0(\boldsymbol{\xi})&=0\ .
\end{aligned}
\end{equation}
with $c_s=\sqrt{g(D_0+A)}$, and $\lambda=\sqrt{3A/4D_0^3}$. 
The exponent $\gamma$ is taken to be equal to $2$ but it can be any number in the interval $[1,3]$. Computing the approximation $w_h$, we can construct the function 
$\eta_h$ by considering the $L_2$-projections of (\ref{eq:sysode1}) in $V_h^{r_1}$, while the velocity field $\bu_h$ can be recovered by solving the linear system (\ref{eq:petsystem0}).

Substitution of the {\em ansatz} (\ref{eq:ansatz}) into the BBM-BBM system (\ref{eq:bbm}) with flat bottom $D(\bx)=D_0$ and following similar computations as before leads to the system
\begin{equation}\label{eq:sysode2}
\begin{aligned}
-c_s\eta+(D_0+\eta)w+\frac{1}{6}c_sD_0^2\nabla_{\xi\zeta}^2\eta&=0\ ,\\
-c_sw+g\eta+\frac{1}{2}w^2+\frac{1}{6}c_sD_0^2w_{\xi\xi}&=0 \ ,\\
\tilde{w} &=0 \ .
\end{aligned}
\end{equation}
In this case we cannot eliminate the variable $\eta$ as in the Peregrine system so instead of solving a scalar ordinary differential equation we practically need to solve a system of equations that involves also the variable $\eta$. The system (\ref{eq:sysode2}) can be written in the form 
\begin{equation}\label{eq:sysode3}
\tilde{\mathcal{L}}\bw=\tilde{\mathcal{N}}(\bw)\ ,
\end{equation}
with $\bw=(\eta, w, \tilde{w})^T$ and
\begin{equation}\label{eq:sysode4}
\tilde{\mathcal{L}}=\begin{pmatrix}
c_s-\frac{1}{6}c_sD_0^2\nabla_{\xi\zeta}^2 & -D_0 & 0\\
-g & c_s-\frac{1}{6}c_sD_0^2\partial^2_{\xi} & 0 \\
0 & 0 & 1
\end{pmatrix}\quad \mbox{and}\quad \tilde{\mathcal{N}}({\bf w})=\begin{pmatrix} \eta w  \\ \frac{1}{2}w^2  \\ 0
\end{pmatrix}\ .
\end{equation}
The Galerkin finite element method for solving system (\ref{eq:sysode3})--(\ref{eq:sysode4}) is then expressed as follows: Seek for the approximation ${\bf w}_h:=(\eta_h,w_h)^T\in U_h^{r_1}\times V_h^{r_2}$, such that
\begin{equation}\label{eq:petvbbm1}
\tilde{\mathcal{L}}_h({\bf w}_h,\boldsymbol{\chi})=(\tilde{\mathcal{N}}({\bf w}_h),\boldsymbol{\chi})\ ,\quad\mbox{for all}\quad \boldsymbol{\chi}\in V_h^{r_1}\times U_h^{r_2}\ ,
\end{equation}
where $\tilde{\mathcal{L}}_h$ is the bilinear form
\begin{equation}
\begin{aligned}
\tilde{\mathcal{L}}_h({\bf w},\boldsymbol{\chi})&= c_s(\eta,\phi)+\frac{1}{6}c_sD_0^2(\nabla_{\xi\zeta}\eta,\nabla_{\xi\zeta}\phi)-D_0(w,\phi)\\
& + c_s(w,\chi)+\frac{1}{6}c_sD_0^2(w_\xi,\chi_\xi)-g(\eta,\chi)+ (\tilde{w},\psi)\ ,
\end{aligned}
\end{equation}
for all ${\bf w}=(\eta,w,\tilde{w})^T\in V_h^{r_1}\times U_h^{r_2}$ and $\boldsymbol{\chi}=(\phi,\chi,\psi)^T\in V_h^{r_1}\times U_h^{r_2}$. Given an initial guess for ${\bf w}_h^0$ for the solitary wave solution of the BBM-BBM system, the Petviashvili method for solving the nonlinear system of equations (\ref{eq:petvbbm1}) is then defined similarly to the analogous method for the Peregrine system:
\begin{equation}
\tilde{\mathcal{L}}_h({\bf w}^{n+1}_h,\boldsymbol{\chi})=\tilde{M}_n^\gamma(\tilde{\mathcal{N}}({\bf w}^{n}_h),\boldsymbol{\chi}), \qquad n=0,1,\cdots\ ,
\end{equation}
for all $\boldsymbol{\chi}\in V_h^{r_1}\times U_h^{r_2}$ and $\tilde{M}_n$ the scalar  
\begin{equation}
\tilde{M}_n=\frac{\tilde{\mathcal{L}}_h({\bf w}^n_h,{\bf w}^n_h)}{(\tilde{\mathcal{N}}({\bf w}^n_h), {\bf w}^n_h)}\ .
\end{equation}
The initial guess for the solitary wave solution can be taken again to be (\ref{eq:initialguesw}). 

As a stopping criterion for the Petviashvili method we consider the norm of the normalized residual to be less than a prescribed tolerance $\delta$. Specifically, we consider the stopping criderion
$$R_n=|\tilde{\mathcal{L}}_h(\bw_h^{n},\bw_h^n)-(\tilde{\mathcal{N}}(\bw^{n}_h),\bw^{n}_h)|/\|\bw_h^n\|_2<\delta\ ,$$ 
with $\delta$ appropriate tolerance.  In this paper we take $\delta=10^{-5}$. In all the experiments we tested the convergence tolerance was satisfied in less than 20 iterations.   
For example,  when we consider the horizontal channel $[-20,30]\times [-1,1]$ and we seek for solitary waves traveling in the direction $\boldsymbol{\alpha}=(1,0)^T$ with depth $D_0=1$ and $g=1$ the Petviashvili method converged to a solitary wave of amplitude $A=0.3$ with tolerance $\delta=10^{-5}$ in $5$ iterations for the BBM-BBM system and in $16$ iterations for the Peregrine system. For this experiment we considered an unstructured grid of $9,510$ triangles with maximum diameter $0.24$ and minimum $0.09$. The propagation of these solitary waves is presented in the next section.

\begin{remark}
The exact relationship between speed and amplitude is known only for the Peregrine system \cite{AD2012}, which is
$$c_s= \frac{\sqrt{6}(D_0+A)}{\sqrt{3D_0+2A}}\frac{\sqrt{gD_0(D_0+A)\log\left(\frac{D_0+A}{D_0}\right)-gD_0A}}{A}\ .$$
For the BBM-BBM system, where there is not any known speed-amplitude relationship,  the appropriate solitary wave can be achieved using a continuation method.
\end{remark}

\subsection{Solitary wave propagation}

\begin{figure}[ht!]
  \centering
  \includegraphics[width=0.47\columnwidth]{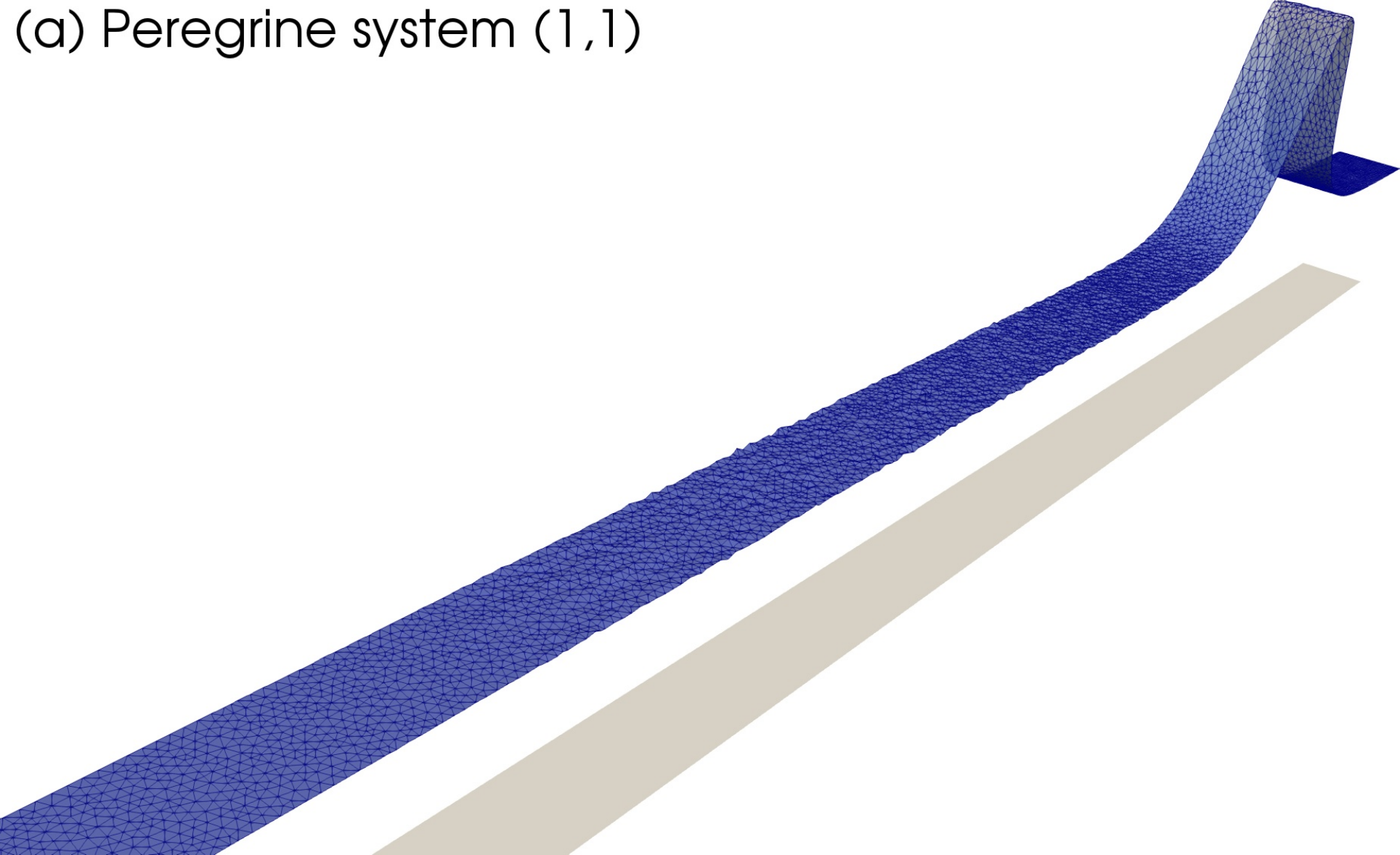}
  \includegraphics[width=0.47\columnwidth]{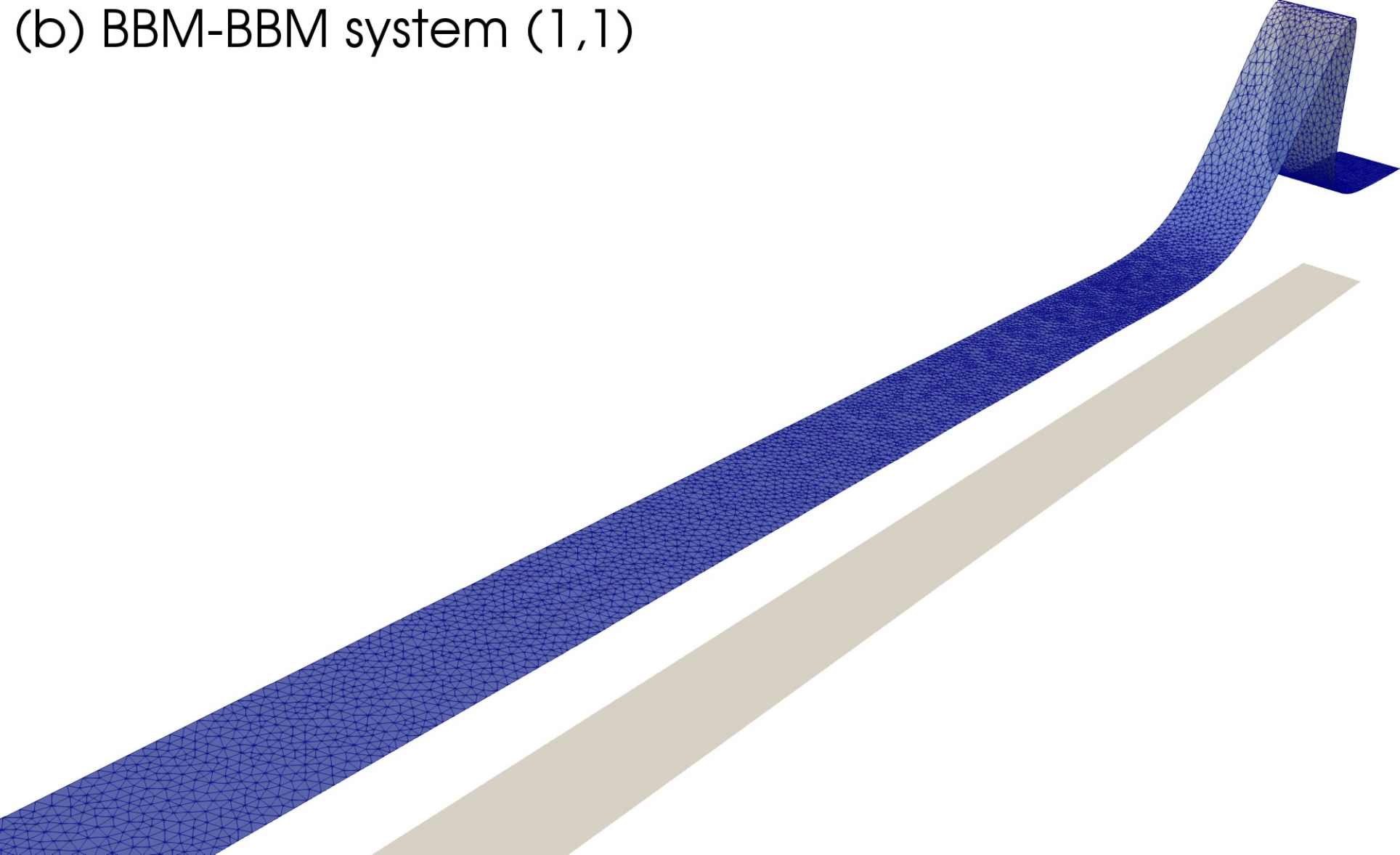}
  \includegraphics[width=0.47\columnwidth]{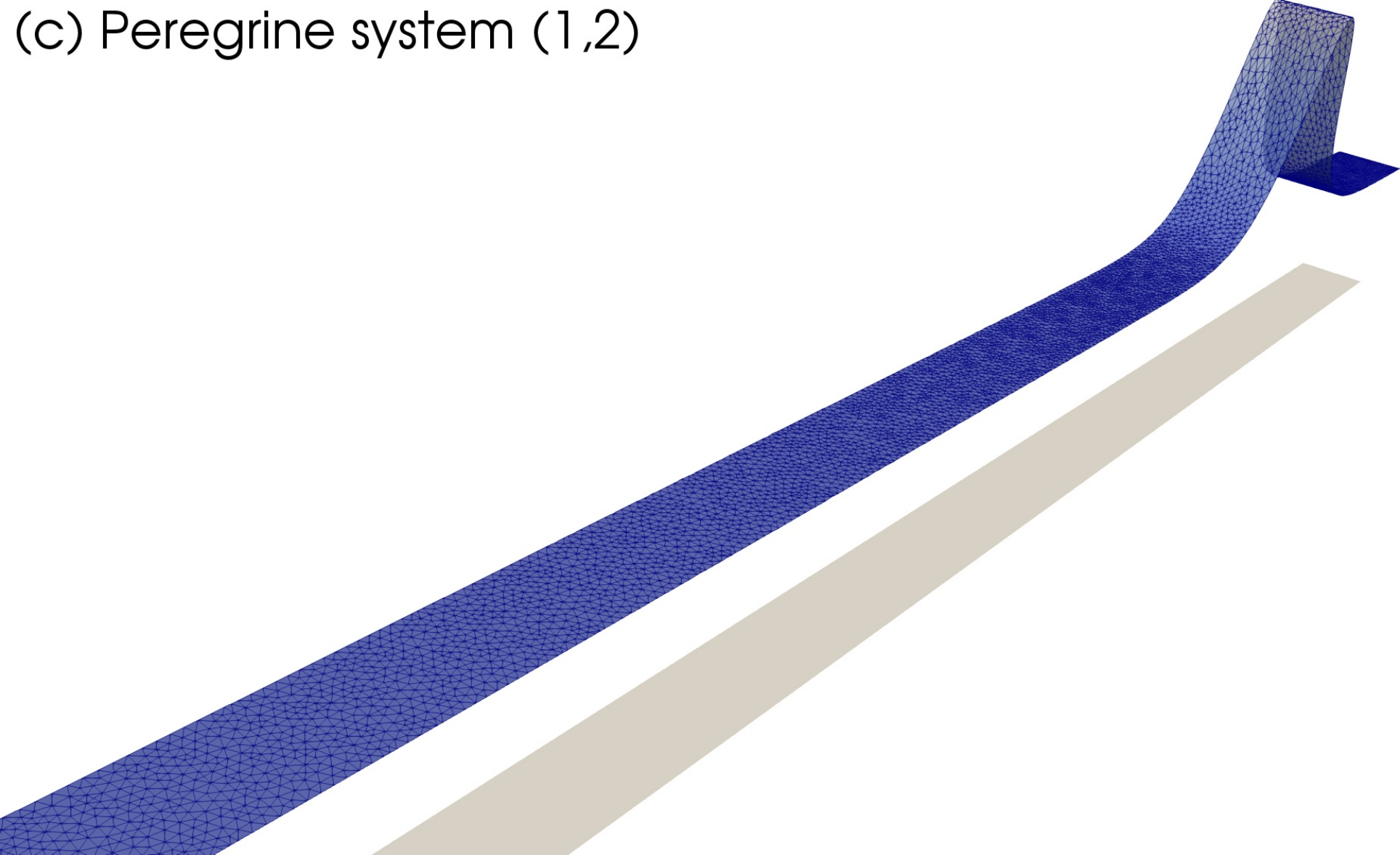}
  \caption{Propagation of a solitary wave of amplitude $A=0.3$ up to $T=25$. ($D=1$, $g=1$)}
  \label{fig:propagation}
\end{figure}

\begin{figure}[ht!]
  \centering
  \includegraphics[width=0.8\columnwidth]{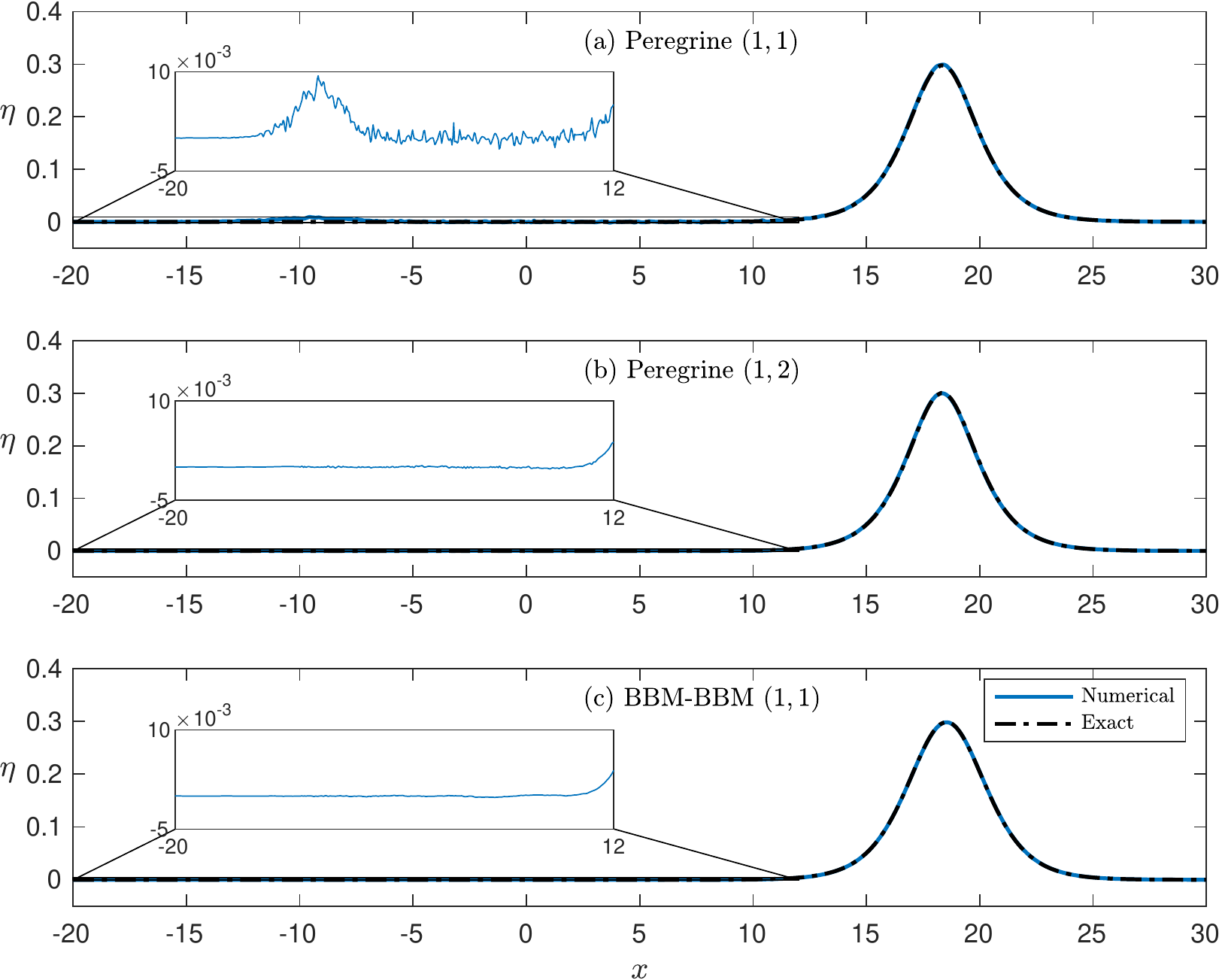}
  \caption{Cross-section along $x$-axis and comparison between the numerical and the exact solution at $T=25$}
  \label{fig:crossection}
\end{figure}

The finite element method has been applied to the Peregrine system before using linear elements in several works. While in \cite{AQ1998,Antunes1993,SSS2004} the propagation of small amplitude solitary waves has been reported to be approximated in a satisfactory way using linear elements, in \cite{ZW2008}  was shown that significant oscillations can be generated behind large amplitude solitary waves during their propagation along a channel. These large oscillations do not affect the stability of the numerical method but they pollute and decrease the resolution of the numerical solution. In \cite{ZW2008}  the use of the streamline diffusion method was proposed to smooth-out the trailing spurious oscillation. However, the streamline diffusion method  is expensive, hard to implement and also adds unwanted damping to the solution. Studying the propagation of a solitary wave of amplitude $A=0.3$ of Peregrine's system in the channel $[-20,30]\times [-1,1]$ with depth $D=1$ up to time $T=25$ we observed  that using  linear elements for all the dependent variables leads to the generation of significant oscillations behind the solitary pulse as shown in Figure \ref{fig:propagation} (a), see also Figure \ref{fig:crossection}. These oscillations are of the exact same nature with those reported in \cite{ZW2008}.  Repeating the same experiment with the BBM-BBM system for a solitary wave of the same amplitude and with the same elements we didn't observe spurious oscillations as the inversion of the elliptic operator in the mass conservation equation regularizes the solution completely. The regularized solution of the BBM-BBM system is presented in Figure \ref{fig:propagation} (b).
Increasing also the degree of the finite elements for the approximation of the velocity field $\bu$ we observe that the spurious oscillations are eliminated. Of course because the initial conditions are not exact traveling waves to the numerical model but just approximations, it is expected that the solitary pulse will be followed by a small amplitude trailing tail which will be of the order degree of the approximation. Thus taking a finer grid or increasing the degree of the finite elements the trailing tails will be decreased in magnitude.

In Figure \ref{fig:crossection} we present a cross-section along the $x$-axis of the three solutions shown in Figure \ref{fig:propagation} at $T=25$. In this figure we can observe clearly  the spurious oscillations generated due to the hyperbolic nature of the mass conservation equation of the Peregrine system with $P_1-P_1$ elements, which are significantly larger than the oscillations generated for the same system using $P_1-P_2$ elements or in the case of the BBM-BBM system with $P_1-P_1$ elements. In the same figure also we observe that in all cases the solitary wave travelled without significant change in shape and speed.
The good behaviour of the finite element method with linear elements for similar BBM-BBM systems has been observed also in \cite{DMS2007,dms2009}. A more interesting observation though here is that the trailing spurious oscillations are eliminated by the use of high-order Lagrange elements. For example by using quadratic elements for the discretization of the velocity field yields a clean solution presented in Figure \ref{fig:propagation} (c). It is noted than in these experiments we considered $\Delta t=0.1$ and the same grid with the one used to generate the solitary waves with the Petviashvili method. 
The previous study of the convergence rates along with the present study of the propagation of the solitary waves suggest that the use of at least quadratic elements for the velocity field of the Peregrine system is mandatory for more accurate and highly-resolved numerical solutions of the Peregrine system. Alternatively, the use of the new regularized BBM-BBM system can be beneficial. In all the following numerical experiments we use linear elements in the approximation of the free-surface elevation and quadratic elements for the approximation of the velocity field.

\section{Numerical experiments}\label{sec:numerexp}

In this section we consider two different experimental settings that demonstrate the ability of the numerical method to handle the propagation and evolution of solitary waves in general two-dimensional domains and variable bottom bathymetric features. Because there are not known exact formulas for the solitary wave solutions propagating without change in shape over a flat bottom $H_0$ of the systems under consideration,  we considered the approximations of the solitary wave obtained using the Petviashvilli method. The experiments we use here, have been used before to validate similar mathematical and numerical models and seem to be appropriate for the regime of weakly nonlinear and weakly dispersive water waves. It is noted that in all the experiments we used $V_h^1$ for the finite element approximation of $\eta$ and $U_h^2$ for $\bf{u}$. 

\subsection{Scattering of solitary wave by a vertical cylinder}

In this first experiment, the propagation of a small amplitude solitary wave and its interaction with a vertical cylinder is numerically investigated. The laboratory experiment was designed and presented in \cite{Antunes1993}. The particular experiment has been used for validation purposes for various Boussinesq systems \cite{KDNS12}.
\begin{figure}[ht!]
  \centering
  \includegraphics[width=0.9\columnwidth]{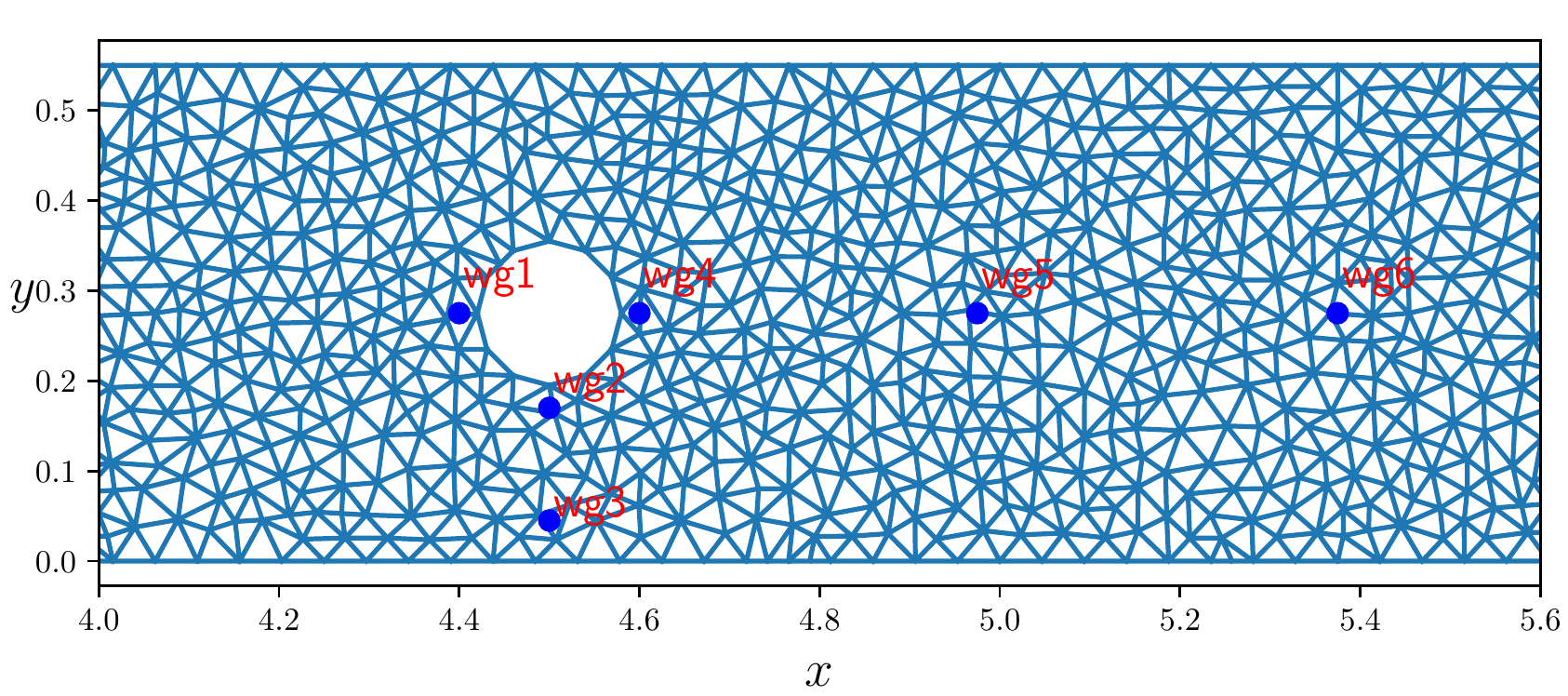}
  \caption{Triangulation around the cylinder and location of the wavegauges}
  \label{fig:triangulation}
\end{figure}
\begin{figure}[ht!]
  \centering
  \includegraphics[width=\columnwidth]{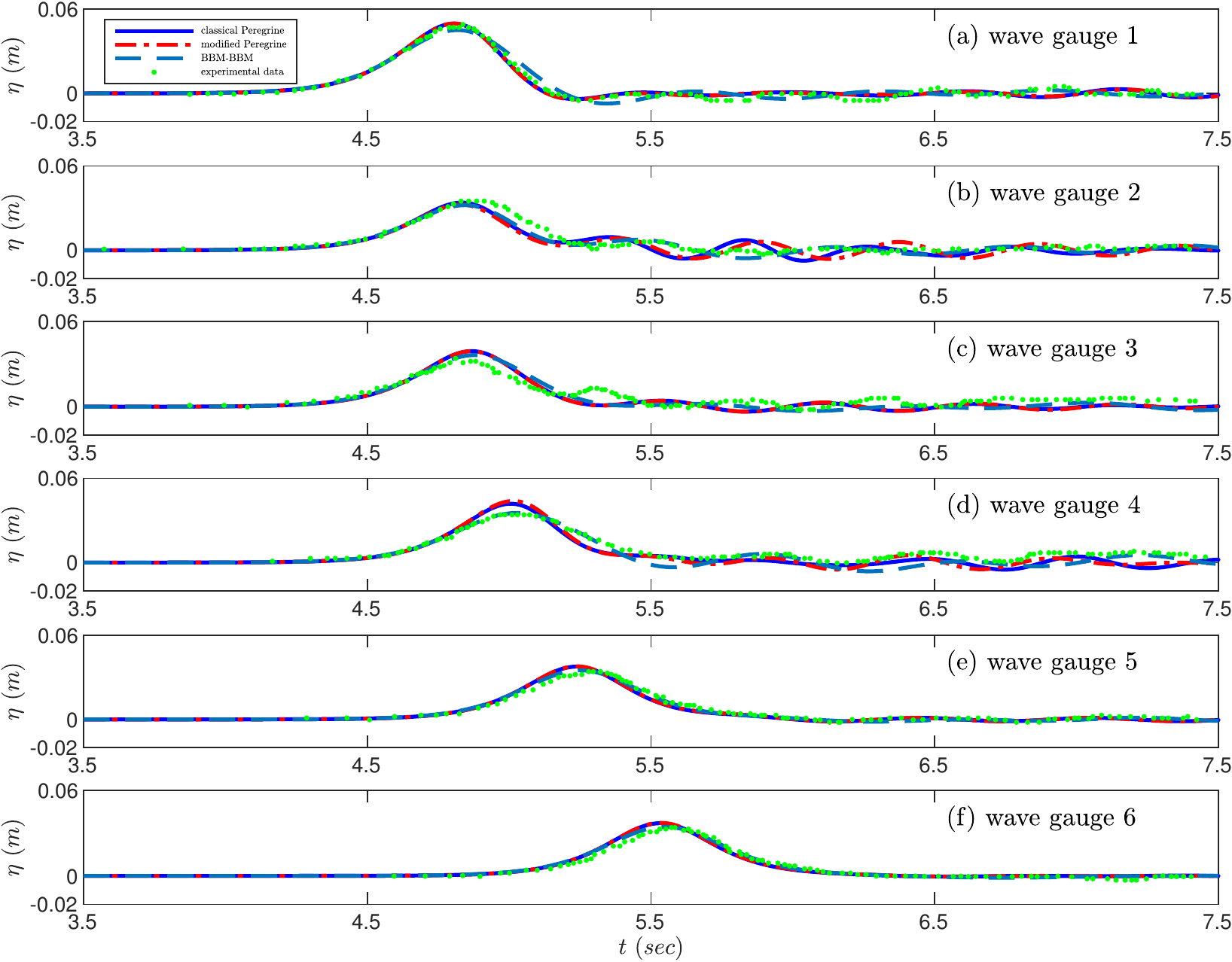}
  \caption{Recorded solution at the six wave-gauges}
  \label{fig:cylinder1}
\end{figure}
\begin{figure}[ht!]
  \centering
  \includegraphics[width=0.41\columnwidth]{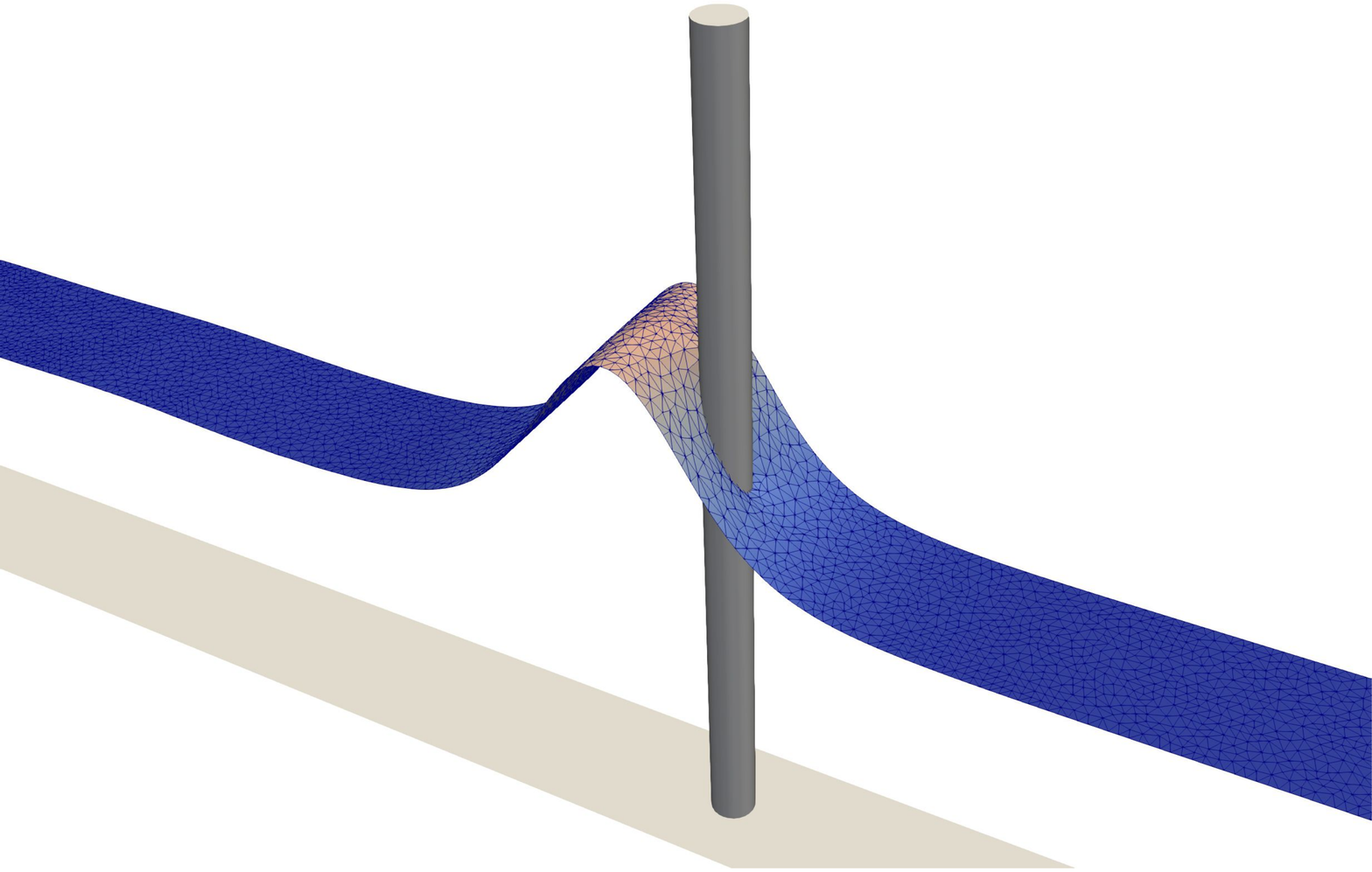}
  \includegraphics[width=0.41\columnwidth]{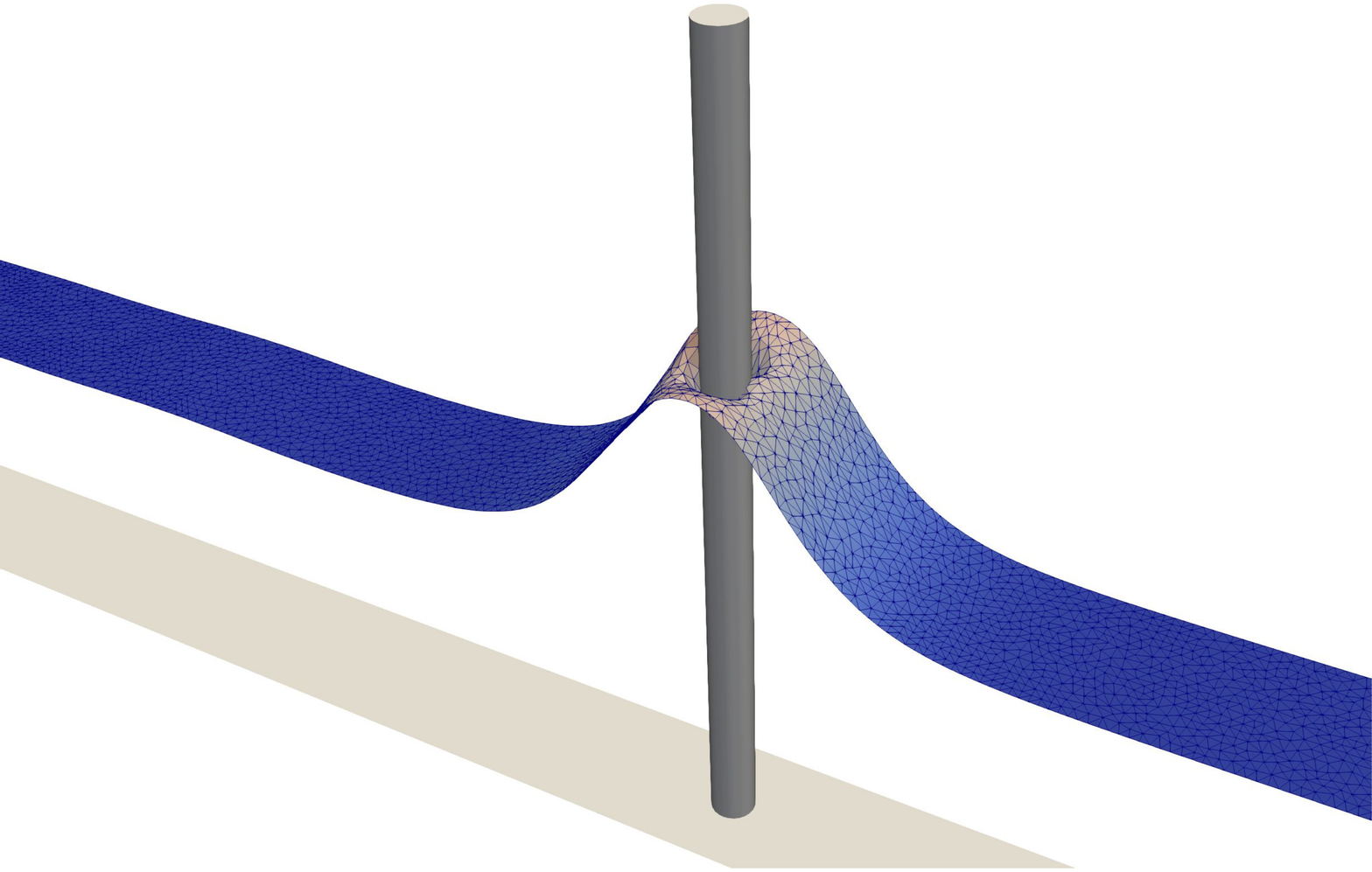}
  \includegraphics[width=0.41\columnwidth]{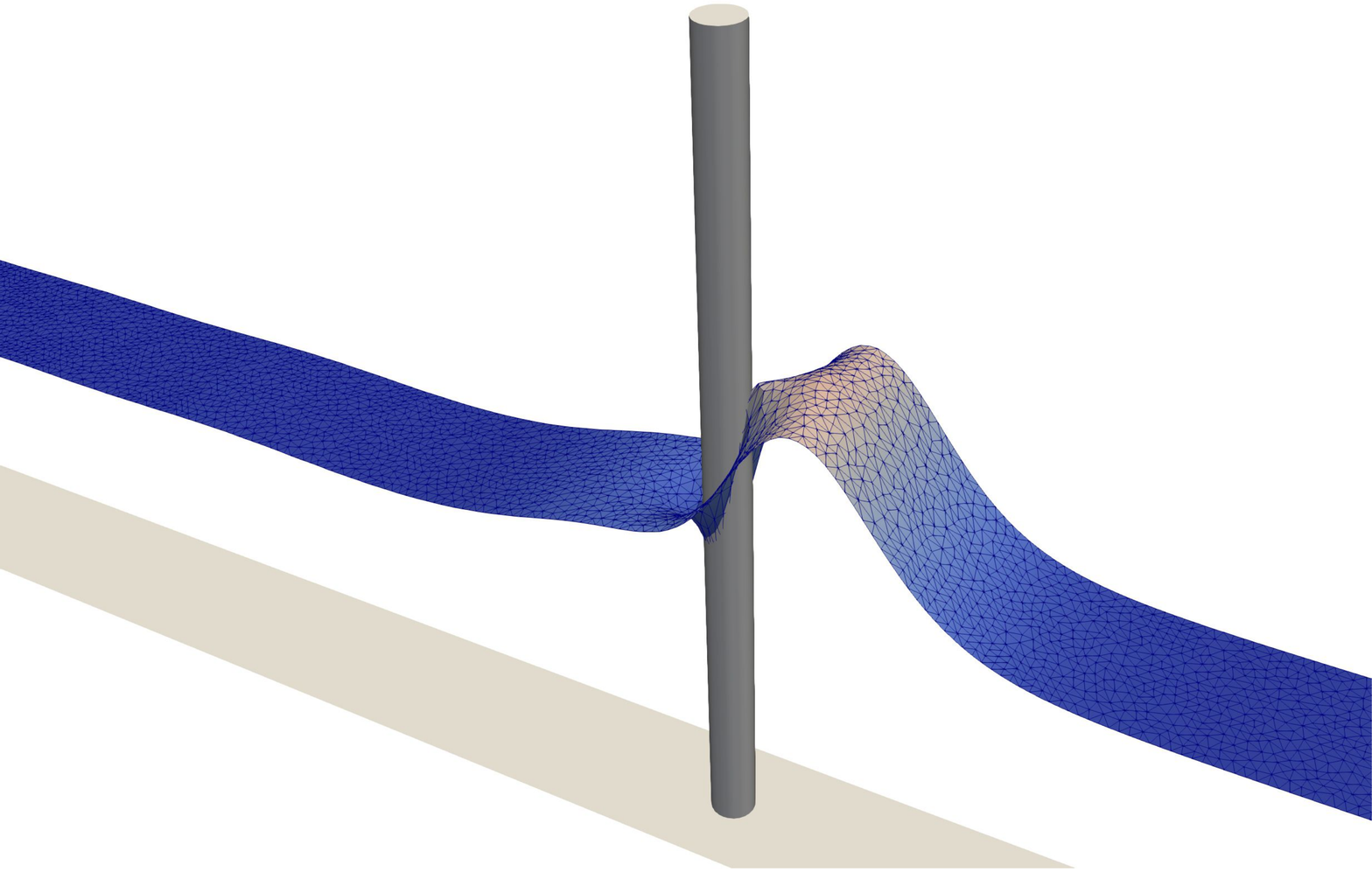}
  \includegraphics[width=0.41\columnwidth]{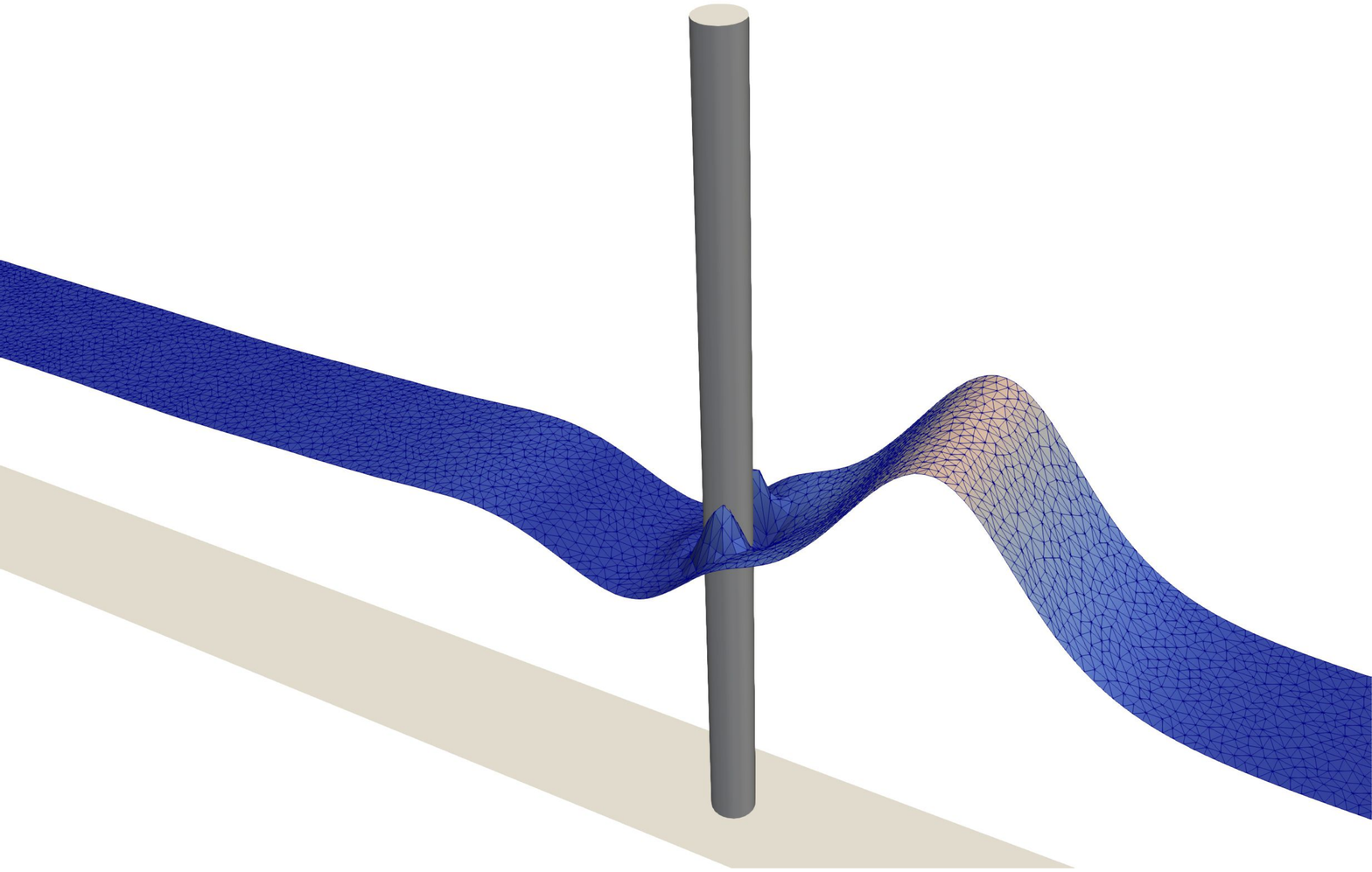}
  \caption{Interaction of the solitary wave with the cylinder at $t=4.688, t=4.928, t=5.152, t=5.472$(left to right, top to bottom), BBM-BBM system}
  \label{fig:cylinder2}
\end{figure}
\begin{figure}[ht!]
  \centering
  \includegraphics[width=0.41\columnwidth]{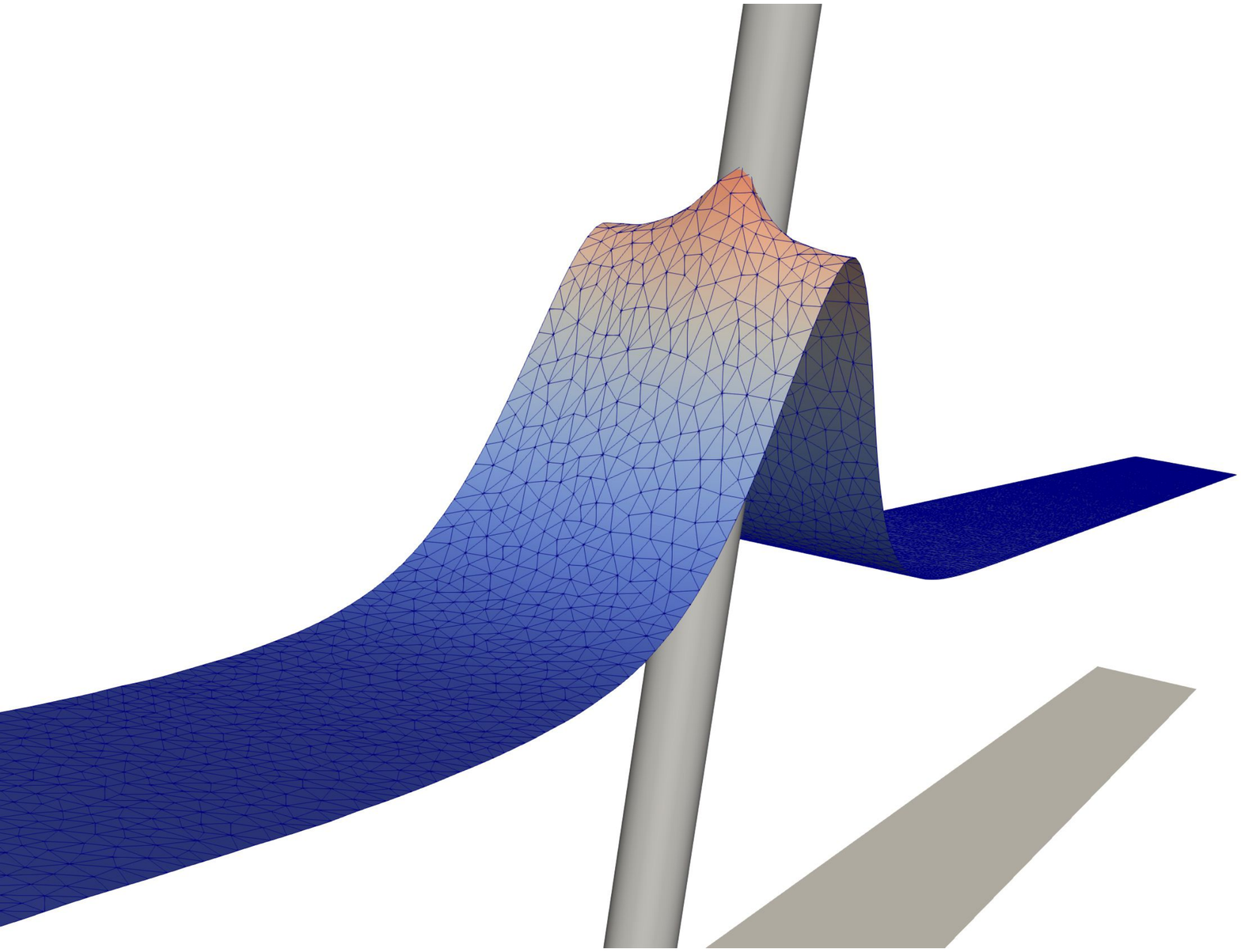}
  \includegraphics[width=0.41\columnwidth]{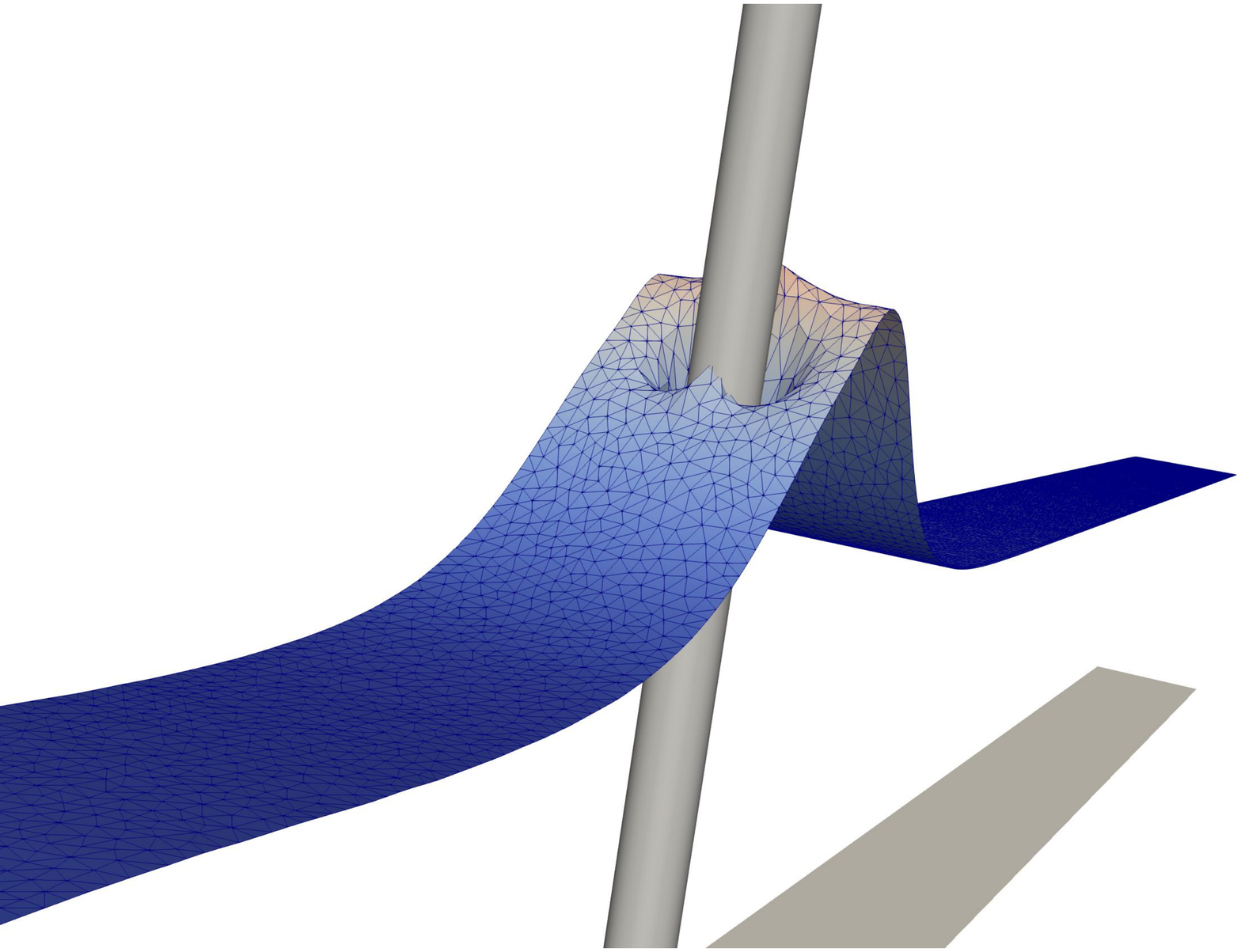}
  \includegraphics[width=0.41\columnwidth]{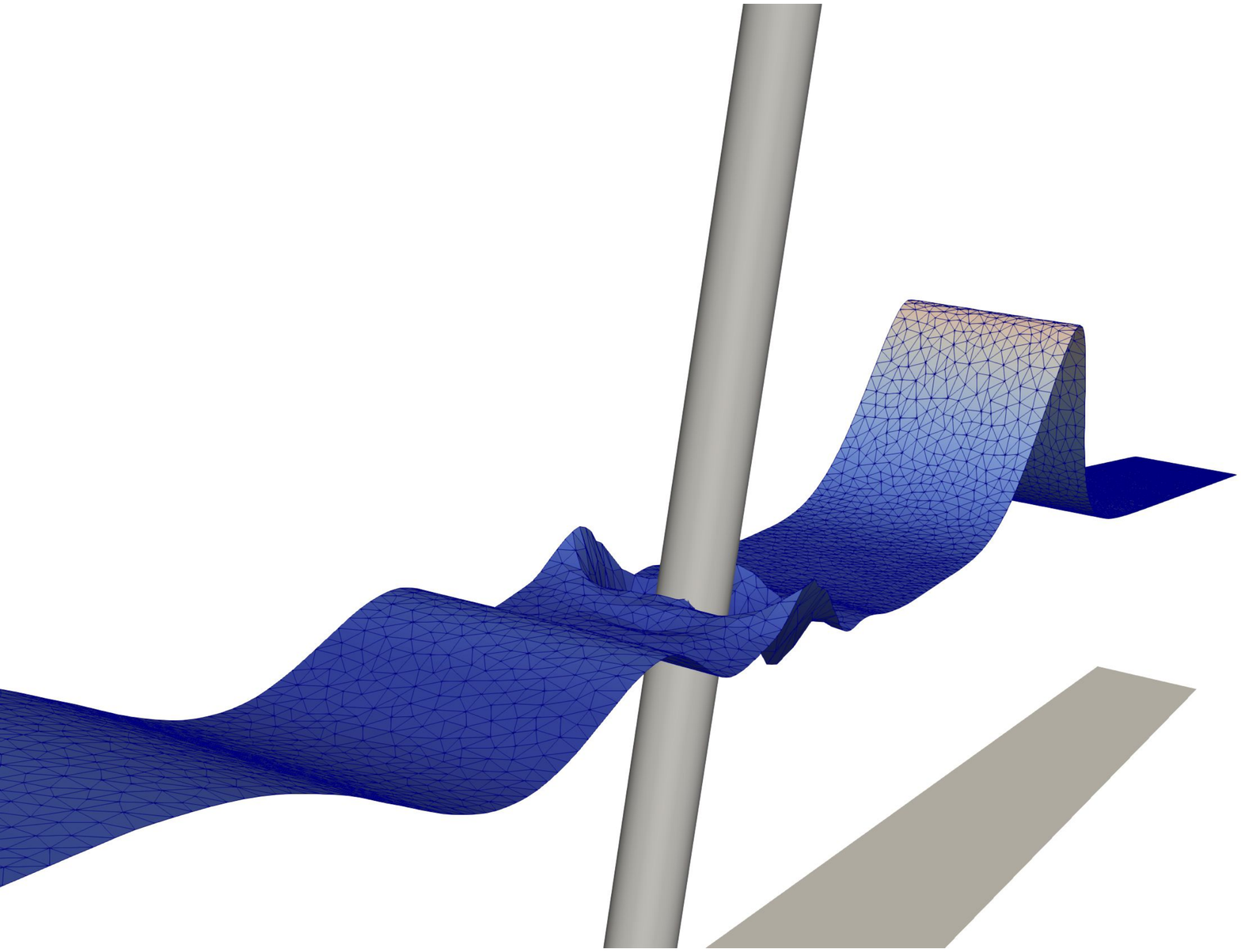}
  \includegraphics[width=0.41\columnwidth]{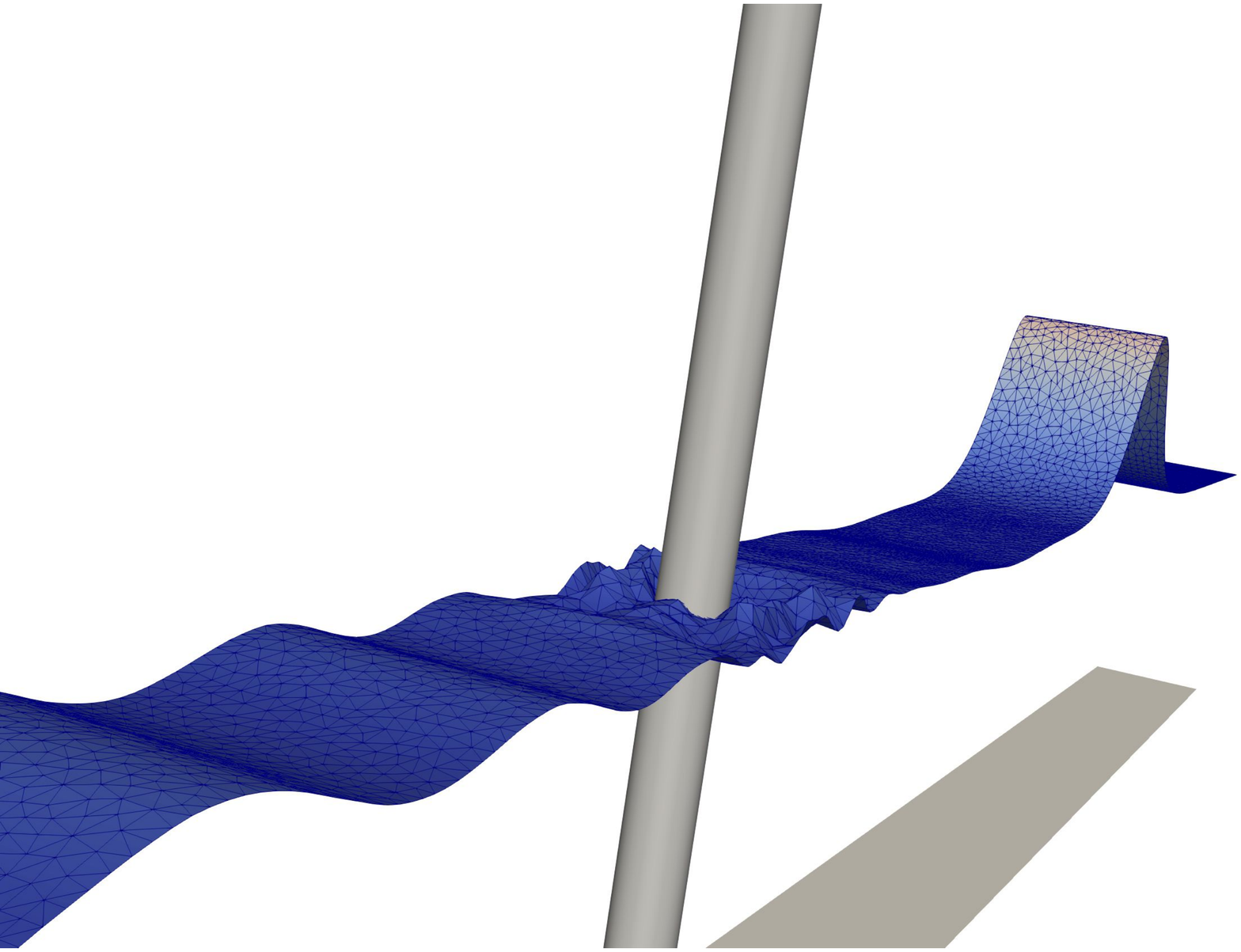}
  \caption{Interaction of the solitary wave with the cylinder at $t=4.832, \ t=5.024, \ t=6.4, \ t=7.536$(left to right, top to bottom), BBM-BBM system}
  \label{fig:cylinder2zoom}
\end{figure}
In this experiment, a rectangular channel of dimensions $[-4,20]\times [0,0.55]$ was used. A vertical cylinder of diameter $0.16~m$ was placed at $(4.5,0.275)$. The solitary wave considered had amplitude $A=0.0375~m$. For the purposes of the numerical experiment we considered the same domain and the initial condition was translated such as its maximum pick was at $2.085~m$. The free-surface elevation $\eta$ was recorded in six wave-gauges located at $\mbox{wg1}=(4.4,0.275)$, $\mbox{wg2}=(4.5,0170)$, $\mbox{wg3}=(4.5,0.045)$, $\mbox{wg4}=(4.6,0.275)$, $\mbox{wg5}=(4.975,0.275)$, $\mbox{wg6}=(5.375,0.275)$. The computational domain was discretised using a regular triangulation consisted of $22,285$ triangles. The timestep is taken $\Delta t=10^{-3}$. The triangulation and locations of the wave-gauges around the cylinder are presented in Figure \ref{fig:triangulation}. 

In every case of Boussinesq system, the solitary wave propagates without significant change in shape (as there are no analytical solutions known for any of the systems it is expected that the approximate solitary wave will propagate as a solitary wave after shedding a small amplitude trailing dispersive tail).  The recorded surface elevation at the wave-gauges agree with the experimental data of \cite{Antunes1993}. The results are presented in Figure \ref{fig:cylinder1}. Apparently, the BBM-BBM system resulted to very accurate results indicating that not only we achieved better convergence of the numerical method but also as a model is a good approximation of the Euler equations for the specific water wave regime. Compared to the numerical results of \cite{Antunes1993} we observe that the results obtained with the current finite element method using triangular grid are better. Possible explanation could be the better resolution of the computational domain and the higher-convergence properties because of the mixed finite element spaces we used.  Compared to the results obtained using the optimised Nwogu's system in \cite{KDNS12} again we observe that in both cases the results are very similar indicating that for the specific experiment, because the solitary wave has relatively small amplitude,  the dispersion relation of the weakly nonlinear structure of the mathematical problems is not significant. 

The interaction of the solitary wave with the cylinder causes the scattering of small amplitude dispersive waves. Due to the small width of the channel, the reflected waves propagate as being like one dimensional waves. On the other hand, as it was observed in similar numerical experiments with wide channel,  these waves propagate as expanding waves in all directions, \cite{DMS2010}.  These features can be seen in a series of snapshots of the free surface $\eta(\bx,t)$ in  Figure \ref{fig:cylinder2}  at various time instances, and with more details in Figure \ref{fig:cylinder2zoom}, where we present the results obtained with the BBM-BBM system.

\subsection{Reflection of shoaling waves}

In the second experiment we consider the shoaling of two solitary waves and their reflection by a vertical wall so as to study the efficiency of the numerical method in experiments with variable bathymetry. The specific experiments have been used before for the validation of various Boussinesq systems \cite{Dodd1998,WB1999,MSM2017}. The propagation of the solitary waves is considered in a channel with dimensions $[-50,20]\times [0,1]$. The bottom topography in this experiment consists of a flat bottom of depth $0.7~m$ for $x<0$ and for all $y$ and a linear bathymetry of slope $1/50$ for $x\in[0,20]$ and for all $y$. The free-surface elevation is recorded by three wave-gauges located at $(0.0,0.5)$, $(16.25,0.5)$ and $(17.75,0.5)$. In the first case, the solitary wave has initial amplitude $0.07~m$. The computational domain was discretised using $14,402$ triangles. The timestep was taken $\Delta t=10^{-3}$.
\begin{figure}[ht!]
  \centering
  \includegraphics[width=0.7\columnwidth]{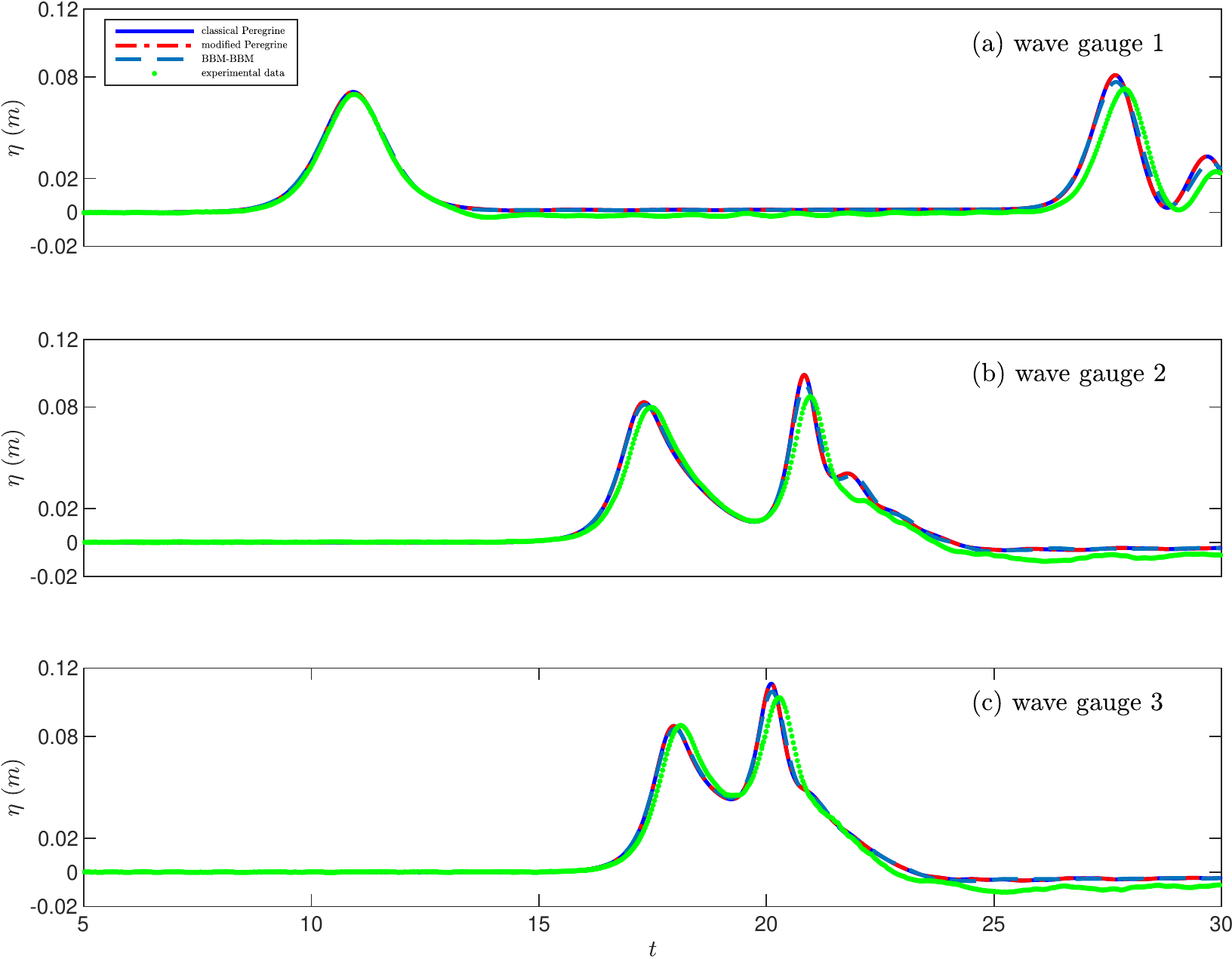}
  \caption{Surface elevation recorded at the three wave-gauges: Case $A=0.07$}
  \label{fig:soal1}
\end{figure}
\begin{figure}[ht!]
  \centering
  \includegraphics[width=0.7\columnwidth]{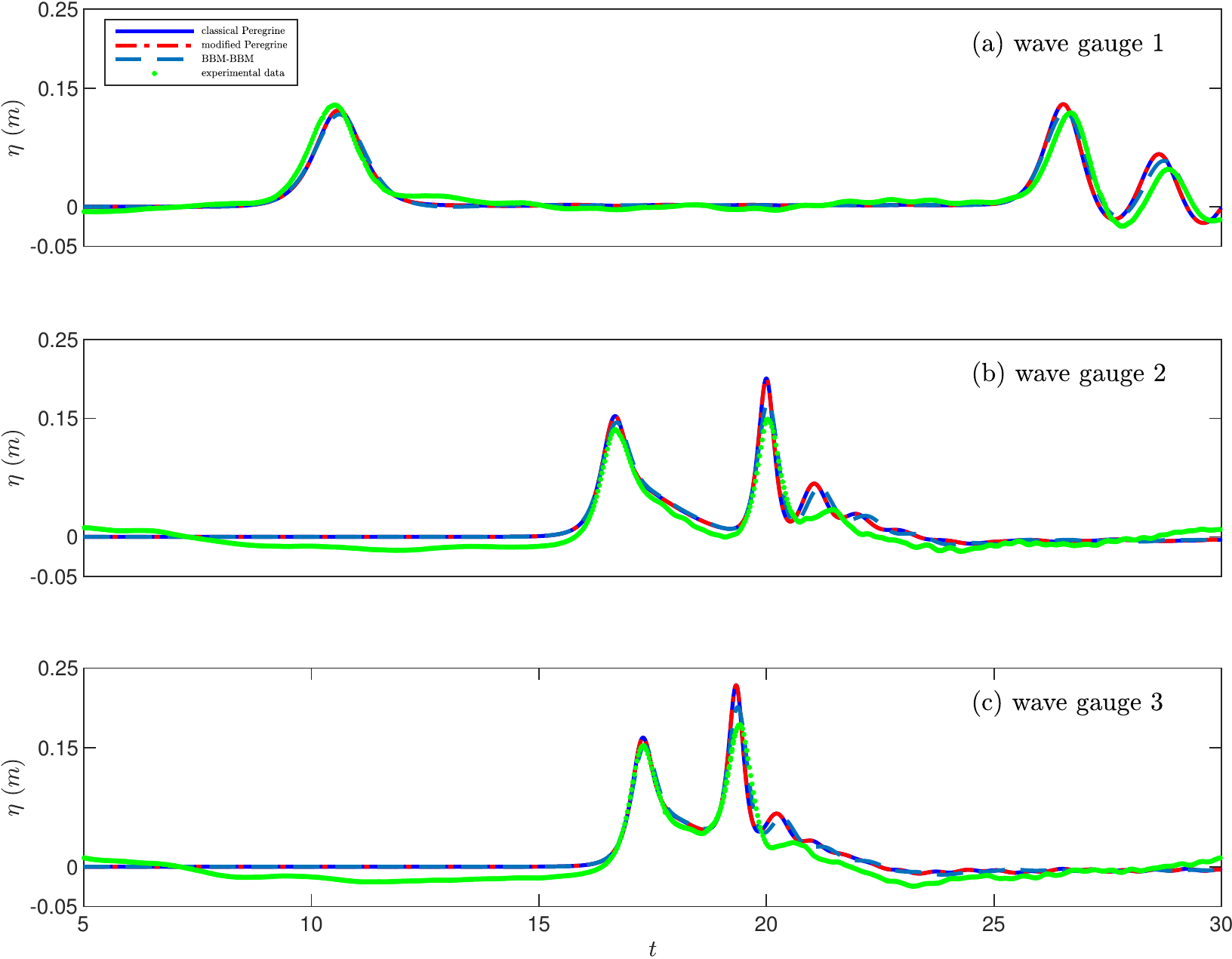}
  \caption{Surface elevation recorded at the three wave-gauges: Case $A=0.12$}
  \label{fig:soal2}
\end{figure}

Because of the small amplitude of the specific solitary wave and the relatively large depth,  the solitary wave propagated without significant change of shape in all cases  and the results fit very well with the laboratory data. Similar observations where made for other Boussinesq systems as well \cite{WB1999,MSM2017} as the specific regime is in favour of all Boussinesq models. It is remarkable that the BBM-BBM system is again closed to the experimental data. The numerical method seems to be efficient during the shoaling and the reflection of the solitary wave. Figure \ref{fig:soal1} presents the recorded surface elevation at the three wave-gauges. 

We also performed the same experiment with a solitary wave of amplitude $0.12~m$ over the same bottom as before. The numerical results obtained at the three wave-gauges are presented in Figure \ref{fig:soal2}. Again the results are very satisfactory especially for the new BBM-BBM system as it appears, even in this case,  the results are slightly better compared to the results obtained using the Peregrine system. 

For the specific experiment the best performance has been observed by the Serre system which is appropriate for experiments with large amplitude solitary wave as there is no small-amplitude assumption in its derivation \cite{MSM2017, KKM2016}.

\subsection{Periodic waves over a submerged bar}

\begin{figure}[ht!]
  \centering
  \includegraphics[width=0.8\columnwidth]{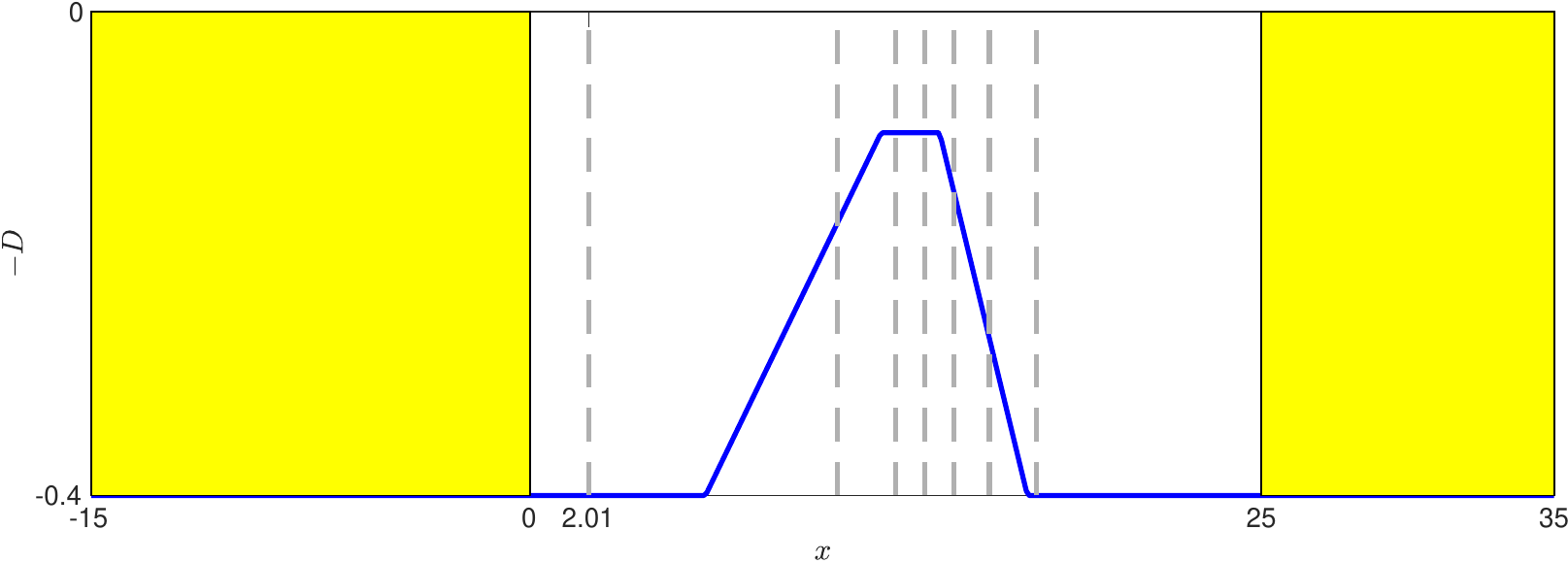}
  \caption{The bathymetry in relation to the sponge layers, wave maker and wave gauges locations}
  \label{fig:bottom}
\end{figure}

In this last experiment we test the numerical method against a laboratory experiment \cite{BB1994} designed to test the nonlinear and dispersive properties of Boussinesq system with variable bathymetry. The specific experiment has been used in several works and is one of the standard benchmarks for numerical models \cite{Dingemans1994,KDNS12,WB1999}. In this experiment periodic waves propagate over a bathymetry defined by the function
$$
D(x,y)=\left\{
\begin{array}{cc}
-0.05x+0.7, & x\in[6,12) \\
0.1, & x\in[12,14) \\
0.1x-1.3, & x\in[14,17] \\
0.4, & \mbox{elsewhere}
\end{array}
 \right.
 $$
and depicted in Figure \ref{fig:bottom}. In this experiment a periodic wave train is generated by a wave-maker. In order to simulate the wave-maker we use the technique suggested in \cite{WB1999} by considering the time-dependent bottom topography incorporated into the Boussinesq systems (\ref{eq:bbmvb}) and (\ref{eq:Peregrinvb}). Specifically, we consider a time-dependent Gaussian bump centered at $x=2.01$ described by the function $\zeta(x,y,t)=ae^{-4(x-2.01)^2}\cos(-\omega t)$ with $\omega=2\pi/k$, $k=2.02$ and $a=0.0095$. The parameters of the moving part of the bottom have been chosen in accordance to \cite{WB1999}. As the bottom oscillates the periodic waves are generated at the free surface with period $k=2.02$. Although the experiment is typically one-dimensional, here we consider the rectangular domain $[-15,35]\times [0,1]$ and a triangulation consisted of $2,500$ triangles. The timestep in this experiment is taken $\Delta t=0.1$. In order to avoid reflection from the lateral boundaries we also apply two sponge layers in the regions with $x\in[-15,0]$  and $x\in[25,35]$ represented by yellow color in Figure \ref{fig:bottom}. The free surface of the water is recorded at  wave-gauges located at $x=10.5$, $12.5$, $13.5$, $14.5$, $15.7$ and $17.3$. The locations of the wave generator and the wave-gauges  are depicted by vertical broken lines in Figure \ref{fig:bottom}.
 \begin{figure}[ht!]
  \centering
  \includegraphics[width=0.8\columnwidth]{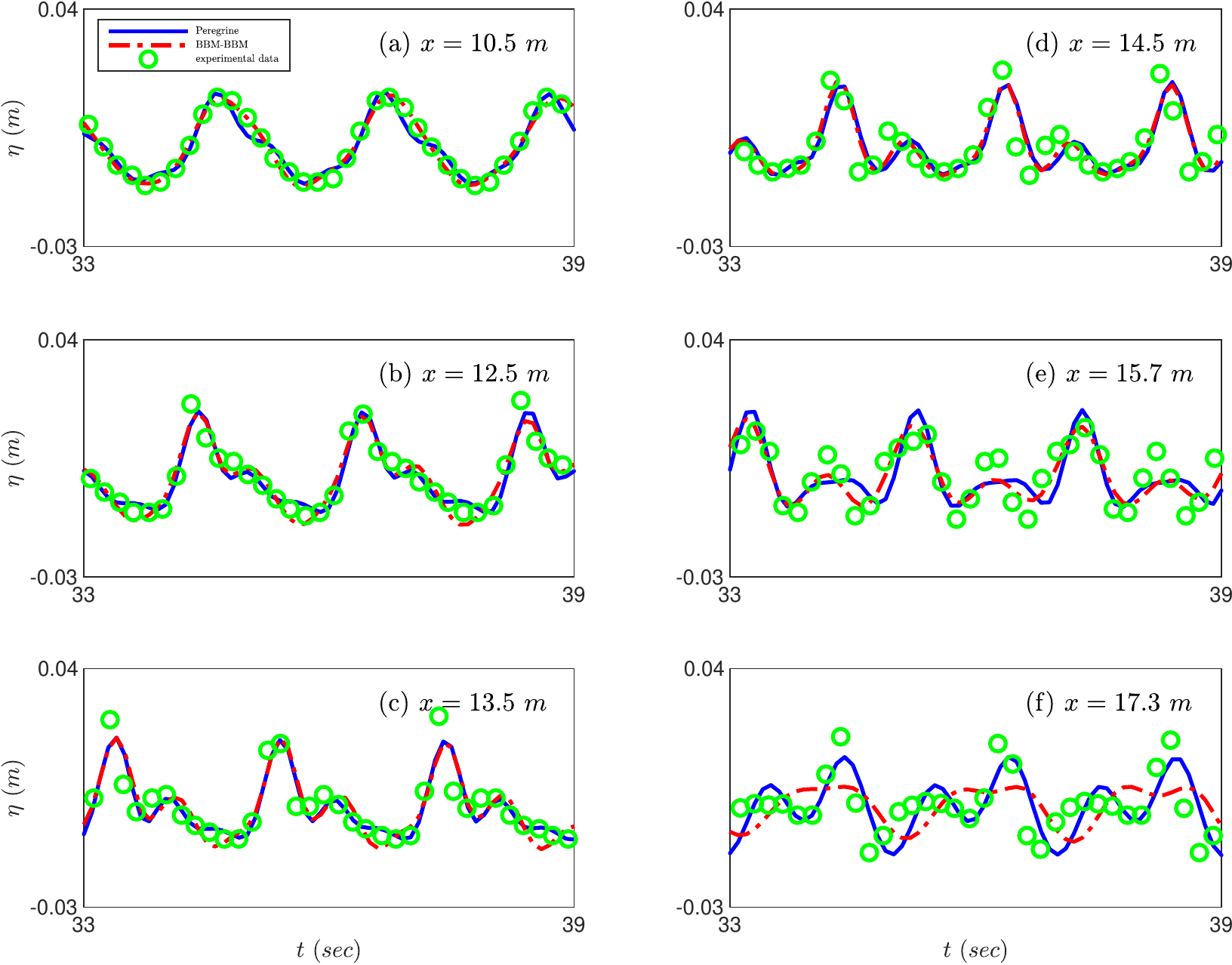}
  \caption{Free-surface elevation recorded at six wave-gauges above the submerged mount}
  \label{fig:pwsb}
\end{figure}
The comparison of the free-surface elevation as predicted by the models and the experimental data is shown in Figure \ref{fig:pwsb}. 
It is known that Boussinesq systems usually fail to predict the correct waves at the wave gauges located further away from the wave maker as the waves generated are nearly breaking and thus outside of the long and small amplitude water wave regime. Boussinesq systems with improved dispersion and nonlinear characteristics, such as Nwogu's system \cite{Nwogu93}, appear to have better behaviour in this experiment although again the errors can be large. In our case both models have almost the same behaviour. As it was also observed in \cite{BB1994} for the Peregrine system, both numerical models start to diverge from the laboratory data after the wave-gauge (d). Moreover, we observe that the numerical results obtained using the BBM-BBM system are in agreement with the numerical results obtained with Peregrine's system until the wave-gauge (e) located at $x=15.7$. After that,  the two models start  to diverge mainly because they have different dispersion relation.
\begin{remark}
The results of the BBM-BBM system can approximate those obtained using Peregrine's system dramatically by considering for example $\theta=\sqrt{1/3}-0.01$, which on one hand has linear dispersion relation very close to the Peregrine system and in addition contains the regularization operator in the mass conservation equation. The numerical solution of the particular system is not presented here since it is very similar to the analogous solution of the Peregrine system.
\end{remark}
\begin{remark}
In all the experiments, the mass conservation error was below $10^{-5}$.
\end{remark}

\section{Conclusions and perspectives}\label{sec:conclusions}

The application of the standard Galerkin finite element method was studied in the case of the Peregrine system with slip-wall boundary conditions. A detailed study of the convergence rates revealed suboptimal convergence properties of the Galerkin method for the particular initial-boundary value problem. Specifically, when we considered the same finite element space for all the dependent variables of the system, then the numerical solution $\eta_h$ converged to the free-surface elevation $\eta$ in the $L_2$-norm with rate of the order $O(h^{r+1/2})$ for elements of degree $r$ with $r$ odd, and $O(h^{r})$ when $r$ was even (when the optimal rate  is $O(h^{r+1})$). The convergence rate of the numerical solution $\bu_h$ to the velocity field $\bu$ was of order $O(h^{r})$ independently of the degree $r$ except the case of linear element $r=1$ where the convergence was optimal. The convergence of both variables $\eta_h$ and $\bu_h$  in $H^1$-norm is one degree smaller than the corresponding rate in $L^2$-norm. Surprisingly the convergence of $\bu_h$ to $\bu$ was optimal in the $\Hdiv$-norm. Mixed finite element methods improved the convergence rates in $\eta$ but not in $\bu$ which remained suboptimal. The suboptimal convergence rates in $\bu$ are due to the significance of the regularization operator which is incomplete-elliptic. This conclusion was drawn after studying alternative systems with different nonlinearity (modified Peregrine system) and by replacing the incomplete-elliptic operator by the standard Laplacian (simplified Peregrine system) using asymptotic arguments.  Although the simplified Peregrine system has better convergence properties, it is hard to be used in situations with slip-wall boundary conditions as it requires the velocity field to be zero on the boundary.

Another problematic situation occurs in the Galerkin method for the Peregrine system when piecewise linear elements are used in two-space dimensions, namely, spurious oscillations can appear behind smooth non-breaking waves propagating over flat bottom. This phenomenon appears because of the hyperbolic structure of the mass conservation equation of the Peregrine system, and can be resolved by using piecewise quadratic (or other high-order) Lagrange elements for the velocity field and/or the free-surface elevation. 

A new Boussinesq system of BBM-BBM type was also proposed as an alternative to Peregrine's system. The momentum conservation equation of the new system has the same structure with the respective equation of the Peregrine system while the mass conservation equation has the classical elliptic regularization operator which regularizes the solution so that no spurious oscillations can be observed even with the use of piecewise linear elements. Moreover, the convergence rates for the free-surface elevation variable are optimal, while the convergence rates for the velocity field is optimal again only in the $\Hdiv$-norm following the same convergence pattern as the Peregrine system.

All the numerical models were tested against standard benchmark experiments within the region of validity of the Boussinesq systems. In particular, we studied the scattering of a line solitary wave by a vertical cylinder, the shoaling of solitary waves and their reflection by a vertical wall and the interaction of periodic waves propagating over an underwater obstacle. The comparison between the new BBM-BBM and the Peregrine systems showed that both systems share the same accuracy against the experimental data for the specific experiments. Thus, the new system can be used in principle as a good alternative to the Peregrine system.

Finally, because there are not any known analytical formulas for solitary wave solutions in closed form of any of the aforementioned systems we employed the Petviashvili method for the computation of numerical approximations of the line solitary waves into the same finite element spaces used for the discretization of the partial differential equations. The convergence and the speed of the Petviashvili method was unexpectedly fast. Further studies for this particular numerical method will follow in the future.

The main conclusions and perspectives are summarised here:
\begin{itemize}
\item The Galerkin finite element method for Boussinesq systems including  other Boussinesq type systems where the mass conservation has the form of a hyperbolic conservation law, requires the use of at least quadratic Lagrange elements in the approximation of the velocity field to achieve satisfactory resolution of the numerical solution.
\item Despite the fact that BBM-BBM type systems have a non satisfactory linear dispersion relationship it seems that in practical problems related to small amplitude long waves, they appear to be as accurate as the Peregrine system.
\item The Petviashvili method applied to the solution of nonlinear equations resulted from the discretization of partial differential equations with finite element methods appear to be very efficient and with good convergence properties.
\end{itemize}

\section*{Acknowledgements}
DM would like to thanks King Abdullah University of Science and Technology for its warm hospitality and support during this project. DM was supported also by the Victoria University of Wellington RSL fund and the Marsden fund administered by the Royal Society of New Zealand.
This work is dedicated to Prof. Vassilios Dougalis on the occasion of his $70^{th}$ birthday. 


\bibliographystyle{plain}
\bibliography{biblio}

\end{document}